\documentclass{article}
\usepackage{amsfonts,latexsym}
\usepackage{amsmath, times}
\usepackage{amssymb}
\usepackage{amsthm}
\usepackage{color}\usepackage[colorlinks=true, citecolor=blue, linkcolor=blue, urlcolor=black]{hyperref}
\usepackage[pagewise]{lineno}
\newcommand{\R}{\mathbb R}

\newcommand{\del}{\partial}
\newcommand{\e}{\varepsilon}
\newtheorem{theorem}{Theorem}[section]
\newtheorem{lemma}{Lemma}

\newtheorem{definition}{Definition}[section]

\newlength{\defbaselineskip}
\newcommand{\setlinespacing}[2]%
{\setlength{\baselineskip}{#1 \defbaselineskip}}

\makeatletter

\@addtoreset{equation}{section} \makeatother 
\begin{document}

\begin{center}
 {\Large   {Existence of sign-changing solutions for a logarithmic weighted Kirchhoff
problem in the whole of $\mathbb{R}^{N}$ with exponential growth non-linearity
 }}
\end{center}
\vspace{0.2cm}
\begin{center}
	  Rached Jaidane $^{(1)}$ \\

	\

	\noindent\footnotesize $^{(1)}$ Department of Mathematics, Faculty of Science of Tunis, University of Tunis El Manar, Tunisia.\\
	Address e-mail: rachedjaidane@gmail.com\\
	
	\
	
\end{center}

\vspace{0.5cm}
\noindent {\bf Abstract.}  In this work, we establish the existence of solutions that change sign at low energy for a non-local weighted Kirchhoff problem in the set $\mathbb{R}^{N}, N>2$. The non-linearity of the equation is assumed to have exponential growth in view of the logarithmic weighted Trudinger-Moser inequalities.
To obtain the existence result, we apply the constrained minimisation in the Nehari set, the quantitative deformation lemma and results from degree theory.\\

\noindent {\bf $2010$ Mathematics Subject classification}:{ $35$J$20$, $35$J$30$, $35$K$57$, $35$J$60$.}

\section{Introduction}
In this paper, we deal with the existence of least energy nodal solutions for the following weighted Kirchhoff problem
\begin{equation}\label{eq:1.1}
  \mathcal{L}_{\omega_{\beta}}(u) =  \displaystyle f(x,u) ~~\mbox{in} ~~\mathbb{R}^{N},   N>2.
\end{equation}
The non-linearity $f(x, t)$ is continuous in $\mathbb{R}^{N}\times \mathbb{R}$ and behaves like $e^{\alpha t^{\frac{N}{(N-1)(1-\beta)}}},~\beta\in[0,1),~~\mbox{as}~t\rightarrow+ \infty$, for some $\alpha >0$ and  $\mathcal{L}_{\omega_{\beta}}$ is the following weighted N-Laplace operator :
$$ \mathcal{L}_{\omega_{\beta}} (u) :=m(\int_{B}(\omega_{\beta}(x)|\nabla u|^{N}dx)\big[-\textmd{div} (\omega_{\beta}(x)|\nabla u|^{N-2}  \nabla u)\big],$$
   where $m$ is a continuous positive function on $(0,+\infty)$, satisfying some conditions which  will be specified later. The weight $\omega_{\beta} (x)$ is given by \begin{equation}\label{wei} \omega_{\beta}(x)= \begin{cases}\left(\log \left(\frac{e}{|x|}\right)\right)^{\beta(N-1)} & \text { if }|x|<1, \\ \chi(|x|) & \text { if }|x| \geq 1,\end{cases}
\end{equation}
where, $0<\beta \leq 1$ and $\chi:[1,+\infty[\rightarrow] 0,+\infty[$ is a continuous function such that $\chi(1)=1$ and $\inf _{t \in[1,+\infty[} \chi(t)>0$. Denote by $E_{\beta}$,  the functional space as follows

$$
E_{\beta}=\left\{u \in L_{\text {rad }}^{\frac{N'}{1-\beta}}\left(\mathbb{R}^{N}\right), \int_{\mathbb{R}^{N}}|\nabla u|^{N} \omega_{\beta}(x) d x<+\infty\right\} .
$$
endowed with the norm
$$
\|u\|=\bigg(\int_{\mathbb{R}^{N}}|\nabla u|^{N} \omega_{\beta}(x) d x\bigg)^{\frac{1}{N}}
$$

For the function $\chi$, we assume that \begin{align}\label{1.1}\int^{+\infty}_{1} \frac{1}{t(\chi(t))^{\frac{1}{N-1}}} d t<+\infty\end{align}In addition, we make the following assumptions: there is a constant  $M>0$ such that\begin{equation} \label{1.2}
    \frac{1}{r^{N^{2}}}\left(\int_{1}^{r} t^{N-1} \chi(t) d t\right)\left(\int_{1}^{r} \frac{t^{N-1}}{(\chi(t))^{\frac{1}{N-1}}} d t\right)^{N-1} \leq M, \forall r \geq 1,
\end{equation}

\begin{equation} \label{1.3}
    \frac{1}{r^{N^{2}}}\left(\int_{1}^{r} t^{N-1} \chi(t) d t\right) \leq M, \forall r \geq 1,
\end{equation}

and

\begin{equation} \label{1.4}
    \frac{\max_{r \leq t \leq 4 r} \chi(t)}{\min_{r \leq t \leq 4 r} \chi(t)} \leq M, \forall r \geq 1.
\end{equation}

The conditions \eqref{1.2},  \eqref{1.3} and  \eqref{1.4}  are
sufficient to guarantee that the weight $\omega_{\beta}$  belongs to the Muckenhoupt’s class $ A_{N}$ ( see\cite{AOJ} ). Here, are some examples of functions $\chi:[1,+\infty[\rightarrow] 0,+\infty[$ satisfying the conditions \eqref{1.2}, \eqref{1.3} and \eqref{1.4}:

$\bullet$ Any function continuous and positive $\chi$ such that $\chi(1)=1$ and

$$
0<\inf _{t \geq 1} \chi(t) \leq \sup _{t \geq 1} \chi(t)<+\infty
$$

$\bullet$ $\chi(t)=t^{\sigma}, 0<\sigma<\frac{N}{N-1}$.

$\bullet$ $\chi(t)=1+\log^{\tau} t,~~\tau>1$.\\

We will now give an overview of the Trudinger-Moser inequalities which correspond to the limiting case of the  Sobolev embedding. Trudinger-Moser type inequalities have been established by S. Pohozaev
\cite{Po}, N. Trudinger \cite{Tru} and V. Yudovich \cite{Ya}. They proved essentially that for some small positive $\alpha>0$, the first order Sobolev space $$ W_{0}^{1,N}(\Omega)=\mbox{closure}\{u\in C_{0}^{\infty}(\Omega)~~|~~\int_{\Omega}|\nabla u|^{^{N}}dx <\infty\}, $$ endowed with the norm \begin{align}\nonumber
\|u\|_{W^{1}_{0}(\Omega)}:=|\nabla u|_{N}=\Big(\int_{\Omega}|\nabla u|^{N}~dx\Big)^{\frac{1}{N}},\end{align} is continuously embedded into the Orlicz space $L_{\phi} (\Omega)$, where $\Omega$ is a smooth bounded domain in $\R^{N}$, $N\geq2$ and $ \phi (t) = e^{\vert t \vert^{\frac{N}{N-1}}}.$

  After, J. Moser\cite{Mos} sharpened this result. More precisely, he showed that for all $ u\in W_{0}^{1,N}(\Omega),$ $$\exp(\alpha \vert u\vert^{\frac{N}{N-1}})\in L^1(\Omega), \quad \forall\, \alpha > 0\: \mbox{and} $$
\begin{equation}\label{eq:Mt}
\sup_{ \|u\|_{W_{0}^{1,N}(\Omega)}\leq 1}
\int_{\Omega}~e^{\alpha|u|^{\frac{N}{N-1} }}dx < C(N)~~~~\Longleftrightarrow~~~~ \alpha\leq \alpha_{N}:=N{\omega}^{\frac{1}{N-1}}_{N-1},
\end{equation}
where $\omega_{N-1}$ is the area of the unit sphere in $\mathbb{R}^{N}$.  The constant $\alpha_{N}$ is sharp in the sense that for $\alpha > \alpha_{N}$
the supremum in \eqref{eq:Mt} is infinite.  Estimates like \eqref{eq:Mt} are now called Trudinger-Moser type inequalities. For related works and applications, we refer the reader to \cite{FMR,FJR,ddr,Liou,Lu} and the references therein.
We mention that great attention has been concentrated on the study of the influence of weights on limiting inequalities of Trudinger-Moser type, we refer to \cite{Adi1, CT} for the effect of power weight in the integral term
on the maximal growth  and \cite{ CR1, CR2, CR3, CRS} for the impact of weights in the Sobolev norm. Kufner \cite{Kuf} established weighted Sobolev spaces and introduced the embedding theory for such weighted Sobolev spaces with  general weight functions. If $\omega \in L^{1}(\Omega)$ be a nonnegative function, we introduce the weighted Sobolev space
\begin{equation} \label{esp}
W_{0}^{1,N}(\Omega,\omega)=\mbox{closure}\{u\in
C_{0}^{\infty}(\Omega)~~|~~\int_{\Omega}\omega(x)|\nabla u|^{^{N}}dx <\infty\}.
\end{equation}

\noindent If the weight function $\omega$ is the logarithmic function, the weighted Sobolev spaces of form \eqref{esp} have a particular significance since they concern limiting situations of such embedding. However, we will need to restrict attention to radial functions to obtain interesting results. So, let us consider the weighted Sobolev space of radial functions.

$$W^{1}_{0,rad}(B,\omega)=closure\{u\in C_{0,rad}^{\infty}(\Omega)~~|~~\int_{\Omega} \omega(x)|\nabla u|^{N}~dx<\infty\}, $$endowed with the norm\begin{equation}\label{norme}
\|u\|_{W^{1}_{0,rad}(B,\omega)}:=|\nabla u|_{N,\omega}=\Big(\int_{\Omega}\omega(x)|\nabla u|^{N}~dx\Big)^{\frac{1}{N}},~~w=\left(\log \left(\frac{e}{|x|}\right)\right)^{\beta(N-1)}.\end{equation}
The first result about the Trudinger-Moser inequalities on Sobolev space with logarithmic weights was
established by Calanchi and Ruf, see \cite{CR1}, where they studied the case when $N = 2$ with Sobolev norm of
logarithmic type. Afterward, they also studied the general case in \cite{CR2}. More precisely, they proved the following result.

\begin{theorem}\cite{CR2} \label{th3}  Let $\beta\in[0,1)$ and let $w$ be given by (\ref{norme}), then
	\begin{equation}\label{eq:1.7}
	\int_{B} e^{|u|^{\gamma}} dx <+\infty, ~~\forall~~u\in W^{1}_{0,rad}(B,\omega),~~
	\mbox{if and only if}~~\gamma\leq \gamma_{N,\beta}=\frac{N}{(N-1)(1-\beta)}=\frac{N'}{1-\beta}
	\end{equation}
	and
	\begin{equation}\label{eq:1.8}
	\sup_{\substack{u\in W^{1}_{0,rad}(B,\omega) \\  |\nabla u|_{N,\omega}\leq 1}}
	\int_{B}~e^{\alpha|u|^{\gamma_{N,\beta} }}dx < +\infty~~~~\Leftrightarrow~~~~ \alpha\leq \alpha_{N,\beta}=N[\omega^{\frac{1}{N-1}}_{N-1}(1-\beta)]^{\frac{1}{1-\beta}},
	\end{equation}
	where $\omega_{N-1}$ is the area of the unit sphere $S^{N-1}$ in $\R^{N}$ and $N'$ is the H$\ddot{o}$lder conjugate of $N$.
\end{theorem}
The existence of solutions to nonlinear weighted elliptic equations involving subcritical and
critical growth of the logarithmic weighted Trudinger-Moser type or weighted Adam's inequalities  has been studied extensively in recent years,
motivated by its applicability in many fields of modern mathematics; see (\cite{ABJ, BJ,CRS, DJ1,DJ2, KJ})  for a survey
on this subject.\\

For unbounded domains, Aouaoui and Jlel have established an extension of Calanchi and Ruff's result to the whole of $\mathbb{R}^{N}$.
More precisely, the authors proved the following weighted Trudinger-Moser inequalities.
\begin{theorem} Let $\beta \in(0,1)$ and $\omega_{\beta}$ be defined by \eqref{wei}. Suppose that the following assumption:
$$(D1)~~ \mbox{There exists a positive constant}~~ C_{0} > 0\mbox{ such that}~~\int^{+\infty}_{1}\left(\int_{r}^{+\infty} \frac{t^{N-1}}{(\chi(t))^{\frac{1}{N-1}}} d t\right)^{N-1}r^{N-1}dr\leq C^{k}_{0},~~\forall ~~k\in \mathbb{N}^{*}$$
~hold, then\\
 For all $\alpha>0$ and $u \in E_{\beta}$, we have

$$
\int_{\mathbb{R}^{2}}\left(e^{\alpha|u|^{\frac{N'}{1-\beta}}}-1\right) d x<+\infty .
$$

Moreover, if $\alpha \leq \alpha_{N,\beta}= N[\omega^{\frac{1}{N-1}}_{N-1}(1-\beta)]^{\frac{1}{1-\beta}}$, then

\begin{equation} \label{1.6}
    \sup _{u \in E_{\beta},\|u\|_{\beta} \leq 1} \int_{\mathbb{R}^{N}}\left(e^{\alpha|u|^{\frac{N'}{1-\beta}}}-1\right) d x<+\infty
\end{equation}
if $\alpha >\alpha_{N,\beta}$, then
\begin{align}
\sup _{u \in E_{\beta},\|u\|_{\beta} \leq 1} \int_{\mathbb{R}^{N}}\left(e^{\alpha|u|^{\frac{N'}{1-\beta}}}-1\right) d x =+\infty
\end{align}
\end{theorem}

We will now give some background on Kirchhoff-type problems. In recent decades, a great many attention has been devoted to the following Kirchhoff-type problem
\begin{equation}\label{eq2}
\displaystyle \left\{
\begin{array}{rclll}
-(a +b \int_{\Omega}|\nabla u|^2)\Delta u&=&  \displaystyle f(x, u)& \mbox{in} & \Omega \\
u&=&0 &\mbox{on }&  \partial \Omega,
\end{array}
\right.
\end{equation}
which is a generalization of a model proposed  by Kirchhoff  \cite{ki}. More precisely, Kirchhoff in 1883 introduced the following equation
\begin{equation}\label{Kir}
	\rho\frac{\partial^2u}{\partial t^2}-\Big(\frac{P_0}{h} +
	\frac{E}{2L}\int_0^L|\frac{\partial u}{\partial
		x}|^2dx\Big)\frac{\partial^2u}{\partial x^2} = 0,
\end{equation}
which extends the classical d’Alembert’s wave equation by considering the effects of the changes in the length of
the strings during the vibrations. The parameters in equation \eqref{Kir} have the following meanings: $\rho$ is the mass density, $P_0$ is the initial tension, $h$ is the area of cross-section, $L$ is the length of the string and  $E$ is the Young modulus of the material. The problem \eqref{Kir} is called nonlocal because the equation contains an integral over $[0, L]$ which makes the study of it interesting.
Note that, \eqref{eq2} can be seen as a stationary version of the following  problem
 \begin{equation}\label{eq:61.4}\left\{
\begin{array}{rclll}
\displaystyle\frac{\partial^{2} u}{\partial t^{2}} -G\bigl(\displaystyle\int_{B}|\nabla u|^{2}dx\bigl)\Delta u  &=& \ f(x,u)& \mbox{in} & B \times (0,T)\\
u &>&0 &\mbox{in}& B\times (0,T)\\
u&=&0 &\mbox{on}&  \partial B \\
u(x, 0)& = &u_{0}(x) &\mbox{in}& B\\
\displaystyle \frac{\partial u}{\partial t}(x,0)&=&u_{1}(x) &\mbox{in}& B,
\end{array}
\right.
\end{equation}
which have focused  the attention of many researchers, mainly as a result of the work of Lions \cite{Li}, where a functional analysis approach was proposed to study it. We mention that non-local problems also arise in other contexts, for example, in biological systems where the function $u$ describes a process that depends on its mean value (e.g., population density), see for example \cite{AC, ACM} and references therein.\\

In (\ref{eq:1.1}), if we set  $m(t)=1$  and $\omega_{\beta}(x)=1$, we find the classical laplacian operateur , which have been extensively studied in the recent years, for example \cite{ABC,A M S,CL}.  Moreover, in the case when $\omega_{\beta}(x) =1,$  and $N=2$, Figueiredo and Severo \cite{FS}  proved the existence of a positive ground state solution to the following Kirchhoff-type problem
$$
 \displaystyle \left\{
      \begin{array}{rclll}
    -m(\int_{\Omega}\vert \nabla u\vert^2 dx)\Delta u &=& \ f(x,u)& \mbox{in} & \Omega \\
u&=&0 &\mbox{on }&  \partial\Omega,
      \end{array}
         \right.$$
where $\Omega$ is a smooth domain of $\mathbb{R}^{2}$, and the nonlinearity $f$ behaves like $\exp\{\alpha t^2\}~\mbox{as}~t\rightarrow+\infty$, for
some $\alpha >0$. The existence result was proved by combining minimax techniques and Trudinger-Moser inequality. A similar result is proved in  \cite{CTW}.

Recently, Abid et al. \cite{ABJ} established the following Kirchhoff-type problem
$$
\left\{
\begin{array}{rclll}
- g(\int_{B} \rho (x)|\nabla u|^{N} dx)\textmd{div} (\rho (x)|\nabla u|^{N-2}  \nabla u) &=& \ f(x,u)& \mbox{in} & B \\
u &>&0 &\mbox{in }& B\\
u&=&0 &\mbox{on }&  \partial B ,
\end{array}
\right.
$$
when $ N\geq 2$, $\rho(x)=(\log\frac{e}{|x|})^{N-1}$ and  the function $f(x, t)$ is continuous in $B \times \mathbb{R}$ and behaves like $\exp\{e^{\alpha{t^{\frac{N}{N-1}}}}\}~~\mbox{as}~~t\rightarrow+\infty$, for some $\alpha>0$. The Kirchhoff function $g$ verifies some conditions.  The authors proved that this problem has a positive ground state solution, using minimax techniques combined with Trudinger-Moser inequality.\\

We now define the notion of subcritical and critical for the nonlinearity $f$ as follows:
\begin{definition}
Let $\gamma:= \gamma_{N,\beta}=\displaystyle\frac{N'}{1-\beta}$,  in view  of inequalities (\ref{eq:1.7}) and (\ref{eq:1.8}), we say that $f$ has subcritical growth at $+\infty$ if
\begin{equation} \label{eq:1.9}
\lim_{|s|\rightarrow +\infty}\frac{|f(x,s)|}{e^{\alpha s^{\gamma}}}=0,~~~~\mbox{for all}~~ \alpha >0
\end{equation}
and $f$ has critical growth at $+\infty$ if there exists some $\alpha_{0}>0,$
\begin{equation}\label{eq:1.10}
\lim_{|s|\rightarrow +\infty}\frac{|f(x,s)|}{e^{\alpha s^{\gamma}}}=0,~~~\forall~\alpha > \alpha_{0}~~~~
\mbox{and}~~~~\lim_{|s|\rightarrow +\infty}\frac{|f(x,s)|}{e^{\alpha s^{\gamma}}}=+\infty,~~\forall~\alpha<\alpha_{0}.\\
\end{equation}
\end{definition}

Concerning sign changing solutions of Kirchhoff type problems, to the best of our knowledge, there are many results in the context, such as \cite{FN, FSJ, GCZ, HY, LDZ, liang4, Sh,  WTC, XTZ, ZYL}. Among them, \cite{S} and \cite{SSLR}, investigated sign-changing solutions of Kirchhoff type problems with  exponential nonlinearities. In particular, Shen in \cite{S} studied the following Kirchhoff type problem

$$
\displaystyle \left\{
\begin{array}{lll}
-(a +b \int_{\R^N}|\nabla u|^N dx)\Delta_N u + V(|x|)|u|^{N-2}u =  \displaystyle f(|x|,u) \: \mbox{in} \: \R^N \\
u \in  W^{1,N}(\R^N),
\end{array}
\right.
$$
with $a, b> 0$, $\Delta_N u=div(\vert \nabla u\vert^{N-2}\nabla u)$, $V(x)$ is a smooth function and $f(|x|,u)$ is critical growth in view of Trudinger-Moser inequalities. By using constraint variational method, the author investigate the existence and asymptotic behavior of least energy sign-changing solutions for the above problem. \\
 Motivated by the above works, in the present paper we investigate the existence of least energy sign-changing solutions for a weighted Kirchhoff problem with subcritical and critical exponential growths at infinity, by using the constraint minimization argument and topological degree theory.

To the best of our knowledge, there are no results concerning the existence of sign-changing solutions for weighted Kirchhoff problems with exponential nonlinearity in the whole of $\mathbb{R}^{N}$.
\section{Assumptions and main results}
Let the space $E_{\beta}$ introduced in \cite{AOJ},
$$
E_{\beta}=\left\{u \in L_{\text {rad }}^{\frac{N'}{1-\beta}}\left(\mathbb{R}^{N}\right), \int_{\mathbb{R}^{N}}|\nabla u|^{N} \omega_{\beta}(x) d x<+\infty\right\},
$$
endowed with the norm
$$
\|u\|=\bigg(\int_{\mathbb{R}^{N}}|\nabla u|^{N} \omega_{\beta}(x) d x\bigg)^{\frac{1}{N}}\cdot
$$

According to Aouiaoui \cite{AOJ}, the embedding $E_{\beta}\hookrightarrow L^{q}(B)$  is continuous for all $q \geq\frac{N'}{1-\beta}$. Moreover,  $E_{\beta}$ is compactly embedded in $L^{q}(B)$  for all $q \geq \frac{N'}{1-\beta} $ \cite{AOJ}. Since $N>2$ and $0<\beta<1$ then these embedding remain valid for all $q\geq N$.

Now, we suppose that the function $m$ is continuously differentiable in $\mathbb{R^{+}}$ and verifies
\begin{description}
	\item[$(M_{1})$]  $m$ is increasing with $m(0)=m_{0}>0$;
	\item[$(M_{2})$] $\displaystyle\frac{m(t)}{t} ~~\mbox{is  nonincreasing for} ~~ t>0.$
\end{description}
From $(M_{1})$ and $(M_{2})$, we can get
	\begin{equation}\label{eq:1.12} M(t+s)\geq M(t)+M(s)~~\forall ~~s,t\geq0, \end{equation}
	where $M(t)=\displaystyle\int^{t}_{0}m(s)ds$,
\begin{equation}\label{eq:1.13}m(t)\leq m(1) + m(1) t,\quad \forall t\geq 0\end{equation}	
and
\begin{equation}\label{eq:1.14} m'(t)t < m(t)\cdot\end{equation}
As a consequence of \eqref{eq:1.14} we have
\begin{equation}\label{eq:1.15}\frac{1}{N}M(t)-\frac{1}{2N}m(t)t~~\mbox{is nondecreasing and positive for}~~t > 0.\end{equation}
		
A typical example of a function $m$ fulfilling the conditions $(M_{1})$ and $(M_{2})$ is given by
$$m(t)=a_{0}+b_{0}t, ~~~~ \forall a_{0}, b_{0} > 0 .$$
Another example is given by $m(t)=1+\ln (1+t)$.\\
Furthermore, we assume the following conditions on the nonlinearity $f(x, t)$.
\begin{description}
\item[$(H_{1})$] $f: \mathbb{R}^{N} \times \mathbb{R}\rightarrow\mathbb{R}$ is $C^1$ and radial in $x$.
  \item[$(H_{2})$] There exist  $\theta > 2N$ such that
  $$ 0 < \theta F(x, t) \leq t f(x,t),\qquad  \forall  (x,t) \in \mathbb{R}^{N} \times \R\setminus \{0\},$$
where $F(x, t) = \displaystyle\int_{0}^{t}f(x,s)ds.$
\item [$(H_{3})$] For each $x\in\mathbb{R}^{N}$,~$\displaystyle t \mapsto \frac{f(x,t)}{|t|^{2N-1}}~~\mbox{is increasing for all}~~t \in \R\setminus \{0\}$.
 \item[$(H_{4})$] $\displaystyle\lim_{t\rightarrow 0}\frac{f(x,t)}{|t|^{N-1}}=0$, uniformly in $x \in\mathbb{R}^{N}$.
\item[$(H_{5})$] There exist $p > 2N$ and $C_p > 1$ such that
$$ sgn(t) f(x,t) \geq C_p \vert t\vert^{p-1}, \quad \mbox{for all} \; (x,t)\in B \times \R,$$
where $sgn(t) = 1$ if $t > 0,$ $sgn(t) = 0$ if $t = 0,$ and $sgn(t) = -1$ if $t < 0.$
\end{description}

\noindent A typical example of a function $f$ satisfying the conditions $(H_{1})$, $(H_{2})$, $(H_{3})$, $(H_{4})$ and $(H_{5})$ is given by
$$f(t)= C_p \vert t\vert^{p-2}t + \vert t\vert^{p-2}te^{t^{\gamma}},\quad \mbox{with} \:p>2N.$$

The energy functional associated with the problem (\ref{eq:1.1}) is defined as follows

 \begin{equation}\label{energy}
\mathcal{E}(u):=\frac{1}{N}M(\|u\|^{N})-\int_{\mathbb{R}^{N}}F(x,u)dx.
\end{equation}
Note that, by the hypothesis ($H_{3}$), for any $\varepsilon>0$, there exists $\delta_{0}>0$
   such that \begin{equation}\label{e1}|f(t)|\leq \varepsilon |t|^{N-1},~~\forall ~~0<|t|\leq \delta_{0}.
   \end{equation}
   Moreover, since $f$ is critical (resp subcritical) at infinity, for every $\varepsilon>0$, there exists $C_{\varepsilon}>0 $ such that
   \begin{equation}\label{e2}~~\forall t\geq C_{\varepsilon}~~|f(t)|\leq \varepsilon \exp( ~\alpha|t|^{\gamma}-1)~~\mbox{with}~~\alpha>\alpha_{0} ~~(\mbox{resp}~~ \forall \alpha>0).
   \end{equation}In particular, we obtain for $q\geq N$,\begin{equation}\label{e3} ~|f(t)|\leq \frac{\varepsilon}{C^{q-1}_{\varepsilon}}|t|^{q-1} \exp(\alpha ~|t|^{\gamma}-1)~~\mbox{with}~~\alpha>\alpha_{0} ~~(\mbox{resp}~~ \forall \alpha>0).
   \end{equation}
 Hence, using   (\ref{e1}), (\ref{e2}), (\ref{e3})  and the continuity of $f$,  for every $\varepsilon>0$, for every $q>N$, there exists  a positive constant $C$  such that \begin{align}\label{imp} |f(t)|\leq \varepsilon |t|^{N-1} +C |t|^{q-1}\big(e^{\alpha ~|t|^{\gamma}}-1\big), ~~~~~~\forall\  t\in\mathbb{R}, ~~\forall~\alpha>\alpha_{0} ~~(\mbox{resp}~~ \forall \alpha>0).\end{align}
 It follows from (\ref{imp}) and $(H_{2})$, that for all $\varepsilon>0$, there exists $C>0$ such that
\begin{equation}\label {eq:1.10}
F(t)\leq \varepsilon|t|^{N}+C |t|^{q}\big(e^{\alpha~|t|^{\gamma}}-1\big),~~~~~~\mbox{for all}~~t,\forall~\alpha>\alpha_{0} ~~(\mbox{resp}~~ \forall \alpha>0).
\end{equation}

 So, by (\ref{1.6}) and  (\ref{eq:1.10}) the functional $\mathcal{E}$ given by (\ref{energy}), is well defined. Moreover, by standard arguments,  $\mathcal{E}\in  C^{1}(E_{\beta},\mathbb{R})$ and  we have
\begin{equation}\label{de1}
\langle \mathcal{E}'(u),v\rangle:=m(\|u\|^{N})\big[\int_{\mathbb{R}^{N}}\big(\omega(x)|\nabla u|^{N-2}~\nabla u~\nabla v\big)~dx\big]-\int_{\mathbb{R}^{N}}f(x,u)~ v~dx. 
\end{equation}
Now, we give the definitions of weak and nodal solutions to problem (\ref{eq:1.1}).
\begin{definition}\label{def1} A function $u$ is called a weak solution to (\ref{eq:1.1}) if $u \in E_{\beta}$ and
\begin{equation*}
m(\|u\|^{N})\big[\int_{\mathbb{R}^{N}}\big(\omega_{\beta}(x)|\nabla u|^{N-2}~\nabla u~\nabla v\big)~dx\big]=\int_{\mathbb{R}^{N}}f(x,u)~ v~dx,~~~~~~\mbox{for all }~~v   \in E_{\beta}.
\end{equation*}
\end{definition}
 We define the sign-changing Nehari set as
$$ \mathcal{N} := \{ u \in E_{\beta}, u^{\pm} \neq 0\; \mbox{and}\; \langle \mathcal{E}'(u),u^+\rangle = \langle \mathcal{E}'(u),u^-\rangle = 0\},$$
where $u^+(x) := \max\{u(x),0\}$ and $u^-(x) := \min\{u(x),0\}$.
We also give the following definition of the so called nodal solutions and least energy sign-changing solution of problem (\ref{eq:1.1}).
\begin{definition}
	$v \in E_{\beta}$ is called a nodal or sign-changing solution of problem (\ref{eq:1.1}) if $v$ is a solution of problem (\ref{eq:1.1}) and $v^\pm \neq 0$ in $\mathbb{R}^{N}.$
	
	\noindent $v \in E_{\beta}$ is called least energy sign-changing solution of problem (\ref{eq:1.1}) if $v$ is a sign-changing solution of (\ref{eq:1.1}) and
	$$\mathcal{E}(v) = \inf\{\mathcal{E}(u) :  \mathcal{E}'(u)=0, \; u^\pm \neq 0 \,\mbox{a.e. in}\, \mathbb{R}^{N}\}.$$
\end{definition}
Our approach is to find sign-changing solutions that minimise the associated energy functional $\mathcal{E}$ from the set of all sign-changing solutions to the problem (\ref{eq:1.1}). Obviously, any solution of the problem (\ref{eq:1.1}) with sign change is in the set $\mathcal{N}$. From \eqref{de1}, one has
\begin{equation*}
\langle \mathcal{E}'(u), u^\pm\rangle =  m(\|u\|^{N})\|u^\pm\|^N - \int_{\mathbb{R}^{N}}f(x,u^\pm)~ u^\pm~dx. 
\end{equation*}
Note that, for every $u = u^+ + u^- \in E_{\beta},$ by  using (\ref{eq:1.12}), it is easy to see that

$$ \mathcal{E}(u) \geq \mathcal{E}(u^+) + \mathcal{E}(u^-).$$
Similarly, from the condition $(M_{1})$, we can obtain
$$\langle \mathcal{E}'(u), u^+\rangle \geq\langle \mathcal{E}'(u^+), u^+\rangle  \quad \mbox{and}\quad \langle \mathcal{E}'(u), u^-\rangle \geq \langle \mathcal{E}'(u^-), u^-\rangle.$$
Also, it's easy to see that $$\|u^{+}+u^{-}\|^{N}=\|u^{+}\|^{N}+\|u^{-}\|^{N}.$$

Our first result establishes the existence of sign-changing solutions for \eqref{eq:1.1} in the exponential subcritical cas

\begin{theorem}\label{th1}
	 Assume that $f(x,t)$ has a subcritical growth at $+\infty$ and satisfies the conditions $(H_{1})$, $(H_{2})$, $(H_{3})$  and $(H_{4})$.  Suppose that the function $m$ satisfies $(M_{1})$  and $(M_{2})$. Then problem (\ref{eq:1.1}) admits a least energy sign-changing solution $v \in \mathcal{N}$ with
	 $\mathcal{E}(v) = \inf_{\mathcal{N}} \mathcal{E}(u)$.
	
 In particular, if $m(t)=a+bt, a, b>0$, then problem (\ref{eq:1.1}) admits a least energy sign-changing solution $v \in \mathcal{N}$ with precisely two nodal domains.

\end{theorem}
The second result proves the existence of a sign-changing solution for \eqref{eq:1.1} in the exponential critical case.

\begin{theorem}\label{th2}
	Assume that $f(x,t)$ has a critical growth at $+\infty$ and satisfies the conditions $(H_{1})$, $(H_{2})$, $(H_{3})$, $(H_{4})$ and $(H_{5})$. Suppose that the function $m$ satisfies $(M_{1})$  and $(M_{2})$. Then there exists $\delta>0$ such that
	\begin{itemize}\item[(i)] problem (\ref{eq:1.1}) admits a least energy sign-changing solution $v \in \mathcal{N}$  provided
 \begin{equation}\label{eq:1.18}
C_p > \max\Big(1, \big(\tau^{\frac{p}{p-N}}\frac{4(p-N)N\theta~~ c_{\mathcal{N}_p}}{m_{0}(\theta-2N)(p-2N)}\bigg(\frac{2(\alpha_0 + \delta)}{\alpha_{N,\beta} }\bigg)^{(N-1)(1-\beta)}\Big)^{\frac{p-N}{N}}\big) >0,
 \end{equation}
 where $C_p$ is the constant given in $(H_5)$, $\displaystyle\tau=\frac{m(1)}{m_{0}}+\frac{m(1)}{m^{2}_{0}}\frac{pN}{p-2N}c_{\mathcal{N}_p}$, $c_{\mathcal{N}_p} = \inf_{\mathcal{N}_p} \mathcal{E}_p(u)>0$,
 $$\mathcal{E}_p(u) :=\frac{1}{N}M(\|u\|^{N}) -\frac{1}{p}\int_{\mathbb{R}^{N}} \vert u\vert^p dx  $$ and
 $$\mathcal{N}_p:= \{ u \in E_{\beta}, u^{\pm} \neq 0\; \mbox{and}\; \langle \mathcal{E}_p'(u),u^+\rangle = \langle \mathcal{E}_p'(u),u^-\rangle = 0\}.$$
 \item[(ii)] if $m(t)=a+bt, a>0 ,b>0$, problem (\ref{eq:1.1}) admits a least energy sign-changing solution $v \in \mathcal{N}$ with precisely two nodal domains provided  $C_{p}$ verifies (\ref{eq:1.18}).

 \end{itemize}
\end{theorem}

This paper is organized as follows. In Section 3, we introduce some  useful Lemmas. Section 4 is devoted to the proof of Theorem \ref{th1}. Finally, in Section 5 we establish some estimates and prove Theorem \ref{th2}.\\

\noindent {\bf Notation:} We shall use the following notations in this paper

\begin{itemize}
\item  $C$  designates a positive constant that can change from one line to the next, and we sometimes index the constant to show how it changes.
\item  $\vert u \vert_p$ denotes the norm in the Lebesgue space $L^p(\mathbb{R}^{N})$ for $p \geq 1$.
\item $\vert u \vert_{p, \omega_{\beta}}$ denotes the norm in the weighted Lebesgue space $L^p(\mathbb{R}^{N}, \omega_{\beta})$ which is defined by
$$|u|_{p,\omega_{\beta}}=\Big(\int_{\mathbb{R}^{N}}\omega_{\beta} (x)|u|^{p}~dx\Big)^{\frac{1}{p}}.$$
\end{itemize}

\section{Some useful and technical lemmas}

In this section, we prove some lemmas which are important to obtain the desired results. To this end, let $u \in E_{\beta}$ with $u^{\pm} \neq 0$, we define the  function
$\mathcal{J}: \mathbb{R}_+\times \mathbb{R}_+\to\mathbb{R}$ and the mapping $\mathcal{I}: \mathbb{R}_+\times
\mathbb{R}_+\to\mathbb{R}^2$, where
\begin{equation}\label{eq:2.1}
\mathcal{J}(s, t)=\mathcal{E}(s u^++t u^-)
\end{equation}
and
\begin{equation}
\label{eq:2.2}
\mathcal{I}(s,t)=\left(\langle \mathcal{E}'(s u^++s u^-),s u^+ \rangle,
 \langle \mathcal{E}'(s u^+ +t u^-),t u^- \rangle\right).
\end{equation}
\begin{lemma}\label{nonempty}
For any $u \in E_{\beta}$ with $u^{\pm} \neq 0$, there is a unique maximum point pair $(s_u, t_u) \in \mathbb{R}_+\times \mathbb{R}_+$ of the function $\mathcal{J}$ such that $s_u u^+ + t_u u^- \in \mathcal{N}.$
\end{lemma}

\begin{proof}The proof of this lemma is a three stage process. The first step is to show that for all $u \in E_{\beta}$ with $u^{\pm} \neq 0$, there exists a pair of positive numbers $(s_u,t_u)$ such that $s_uu^++ t_uu^- \in \mathcal{N}$.	Note that
	\begin{equation}\label{eq:2.3}
	\begin{aligned}
	&\langle \mathcal{E}'(s u^++t u^-),s u^+\rangle
	&=m(\|s u^{+} +tu^{-}\|^{N})\|su^{+}\|^{N}- \int_{\mathbb{R}^{N}}f(x,s u^+)~ su^+~dx
	\end{aligned}
	\end{equation}
	and
		\begin{equation}\label{eq:2.4}
		\begin{aligned}
	&\langle \mathcal{E}'(s u^++t u^-),t u^-\rangle
	&=m(\|s u^{+} +tu^{-}\|^{N})\|tu^{-}\|^{N}- \int_{\mathbb{R}^{N}}f(x,t u^-)~ tu^-~dx.
	\end{aligned}
	\end{equation}

We begin with the critical case. From (\ref{imp}), for every $\varepsilon>0$, there exist constants $C>0$, $\alpha_{0}>0$  such that  for all $q>N$

\begin{equation}\label{eq:2.7}f(x,t)t\leq  \e
	|t |^N+ C|t|^{q} exp(\alpha |t|^{\gamma}-1), \quad \mbox{for all}\; \alpha > \alpha_{0},  q>N \;\mbox{and}\; t\in \R.	\end{equation}
So, we get
	$$\int_{\mathbb{R}^{N}}f(x,s u^+)~ su^+dx \leq \e
	\int_{\mathbb{R}^{N}}|s u^+|^N\,dx + C \int_{\mathbb{R}^{N}}|su^+|^{q} exp(\alpha (su^+)^{\gamma}-1)dx, \quad \mbox{for all}\; \alpha > \alpha_0, \; q>N.$$
Now, from \eqref{eq:2.7}, \eqref{eq:1.7}, $(M_{1})$, Sobolev embedding theorem and H\^older inequality, we have for $s > 0$ small enough  satisfying $s\leq \frac{\alpha_{N,\beta}^{\frac{1}{\gamma}}}{(2\alpha)^{\frac{1}{\gamma}}\|u^{+}\|}.$
		\begin{equation}\label{eq:2.8}
	\begin{aligned}
	&\langle \mathcal{E}'(s u^++t u^-),s u^+\rangle
	\\&\geq m(\|s u^{+} +tu^{-}\|^{N})\|su^{+}\|^{N}-\varepsilon s^{N}
	\int_{\mathbb{R}^{N}}| u^+|^N\,dx - C_\e s^q\int_{\mathbb{R}^{N}}|u^+|^{q} exp(\alpha (su^+)^{\gamma}-1)dx
	\\&\geq m_{0}\|su^{+}\|^{N}- \e s^N
	\int_{\mathbb{R}^{N}}| u^+|^N\,dx\\&\quad - C s^q\big(\int_{\mathbb{R}^{N}}|u^+|^{2q}dx\big)^{\frac{1}{2}} \big(\int_{\mathbb{R}^{N}} exp(2 \alpha s^{\gamma} \|u^+\|^{\gamma} (\frac{u^+ }{\Vert u^+\Vert})^{\gamma}-1)dx\big)^{\frac{1}{2}}
	\\&\geq m_{0} s^N \Vert u^+\Vert^N -\e s^N C_{1} \|u^{+}\|^{N} - C_2 s^q \Vert u^+\Vert^q.
    \end{aligned}
	\end{equation}
	Choose $\e>0$ small enough such that $(m_{0}-\e C_1)>0$. Since $q>N$,  we have that
	$$\langle \mathcal{E}'(s u^++t u^-),s u^+\rangle>0,\quad \mbox{for}\; s\; \mbox{small enough and all}\; t\ge0.$$
	 Similarly, according to \eqref{eq:2.4} and \eqref{eq:2.7}, we get $$\langle \mathcal{E}'(s u^++t u^-),t u^-\rangle>0, \quad \mbox{for}\; t\; \mbox{small enough and all}\; s\ge0.$$
	Hence, there exists $r>0$ such that
	\begin{equation}\label{eq:2.9}
	\langle \mathcal{E}'(ru^++ tu^-),ru^+\rangle>0\quad \text{and} \quad
	\langle \mathcal{E}'(s u^++ru^-),ru^-\rangle>0, \quad \mbox{for all}\; s, t\ge0.
	\end{equation}
	
	On the other hand, by $(H_2)$ and  $(H_4)$ , we can find positive constants $C_3$ and $C_4$  such that
	\begin{align}\label{eq:2.10}
	f(x,t)t \geq C_{3} \vert t\vert^\theta - C_{4}|t|^{N}, \quad \forall (x,t)  \in (\mathbb{R}^{N},\mathbb{R}\setminus \{0\}).
	\end{align}
	Thus, using (\ref{eq:1.13}) we get
	\begin{equation}\label{eq:2.11}
	\begin{aligned}
	&\langle \mathcal{E}'(s u^++t u^-),s u^+\rangle
	\\& =m(\|s u^{+} +tu^{-}\|^{N})\|su^{+}\|^{N}- \int_{\mathbb{R}^{N}}f(x,s u^+)~ su^+~dx
	\\&\leq  m(1) \|su^{+}\|^{N} +m(1)\|s u^{+} +tu^{-}\|^{N}\|su^{+}\|^{N} - C_3 s^\theta
	\int_{\mathbb{R}^{N}}| u^+|^\theta +  C_4 |s|^{N}\int_{\mathbb{R}^{N}}| u^+|^N.
	\end{aligned}
	\end{equation}
	
	Since $\theta >2N$, there exists $R>r$ large enough such that
	\begin{equation}\label{eq:2.12}
	\langle \mathcal{E}'(Ru^++ tu^-),Ru^+\rangle< 0\quad \text{and} \quad
	\langle \mathcal{E}'(s u^++Ru^-),Ru^-\rangle< 0, \quad \mbox{for all}\; s, t \in [r,R].
	\end{equation}
	
	In accordance with Miranda's theorem \cite{mi}, as well as \eqref{eq:2.9}
	and \eqref{eq:2.12}, we can conclude that there exists
	$(s_u,t_u)\in \mathbb{R}_+ \times
	\mathbb{R}_+$ such that $\mathcal{I}(s_u,t_u)=(0,0)$, i.e.,
	$s_uu^++t_uu^-\in \mathcal{N}$.
	
For the subcritical case, \eqref{eq:2.4} and \eqref{eq:2.7} hold for all $\alpha > 0$ and the rest of the proof is the same as in the critical case.\\

	In the second step, we will show the uniqueness of the pair $(s_u,t_u)$. First, we assume that $u = u^++ u^- \in \mathcal{N}$. From Claim 1, we know  that there exists at least  one  positive pair $(s_0, t_0)$ satisfying $s_0u^++t_0u^-\in \mathcal{N}$. Now, we show  that $(s_0,t_0)=(1,1)$ is the unique pair of numbers. Without
	loss of generality, let us assume that $s_0\leq t_0$. It follows from
	\eqref{eq:2.11} and $(M_1)$ that
	\begin{equation} \label{eq:2.13}
	\begin{aligned}
m(s_0^N \Vert u\Vert^N)s_0^N \Vert u^+\Vert^N \leq \int_{\mathbb{R}^{N}}f(x,s_0 u^+)~ s_0u^+~dx.
	\end{aligned}
	\end{equation}
So,
\begin{equation} \label{eq:2.14}
\begin{aligned}
\frac{m(s_0^N\Vert u\Vert^N)}{s_0^N\Vert u\Vert^N}  \leq \frac{1}{\Vert u\Vert^N\Vert u^+\Vert^N} \int_{\mathbb{R}^{N}}\frac{f(x,s_0 u^+)}{(s_0 u^+)^{2N-1}} (s_{0}u^+)^{2N}~dx.
\end{aligned}
\end{equation}	
On the other hand, since $u \in \mathcal{N}$, we have
\[
\langle \mathcal{E}'(u),u^+\rangle=0\quad \text{and} \quad
\langle \mathcal{E}'(u),u^-\rangle=0,
\]
that is,
\begin{equation} \label{eq:2.15}
\frac{m(\Vert u\Vert^N)}{\Vert u\Vert^N} = \frac{1}{\Vert u\Vert^N\Vert u^+\Vert^N}\int_{\mathbb{R}^{N}}\frac{f(x, u^+)}{(u^+)^{2N-1}} (u^+)^{2N}~dx
\end{equation}
and
\begin{equation} \label{eq:2.16}
\frac{m(\Vert u\Vert^N)}{\Vert u\Vert^N} = \frac{1}{\Vert u\Vert^N\Vert u^-\Vert^N}\int_{\mathbb{R}^{N}}\frac{f(x, u^-)}{(u^-)^{2N-1}} (u^-)^{2N}~dx.
\end{equation}
	If $s_0<1$, then from \eqref{eq:2.14}, \eqref{eq:2.15},  $(M_2)$ and
	$(A_3)$, we have
	\begin{equation} \label{eq:2.17}
	\begin{aligned}
	0& <\frac{m(s_0^N\Vert u\Vert^N)}{s_0^N\Vert u\Vert^N} - \frac{m(\Vert u\Vert^N)}{\Vert u\Vert^N}
	\leq \frac{1}{\Vert u\Vert^N\Vert u^+\Vert^N}\int_{\mathbb{R}^{N}}\Big(\frac{f(x,s_0 u^+)}{(s_0 u^+)^{2N-1}}-\frac{f(x, u^+)}{(u^+)^{2N-1}}\Big)(u^+)^{2N}\,dx  \leq 0
	\end{aligned}
	\end{equation}
	which is a contradiction. Hence, $1\leq s_0\leq t_0$.
	
	Arguing similarly by using the equations $\langle \mathcal{E}'(s u^++t u^-),t u^-\rangle =0$ and $\langle \mathcal{E}'(u),u^-\rangle=0$, we obtain that $s_0\leq t_0\leq 1$, which
	implies that $s_0=t_0=1$ and the proof is complete.
	
For the general case, we suppose that $u \notin \mathcal{N}$. Assume there exist two other pairs of positive numbers $(s_1,t_1)$ and
	$(s_2,t_2)$ such that
	\[
	\sigma_1=s_1u^+ +t_1u^-\in\mathcal{N}\quad \text{and}
	\quad \sigma_2=s_2u^++t_2u^-\in\mathcal{N}.
	\]
	Then
	\[
	\sigma_2=\big(\frac{s_2}{s_1}\big)s_1u^+
	+\big(\frac{t_2}{t_1}\big)t_1u^-
	=\big(\frac{s_2}{s_1}\big)\sigma_1^+
	+\big(\frac{t_2}{t_1}\big)\sigma_1^-
	\in\mathcal{N}.
	\]
	Since $\sigma_1\in\mathcal{N}$, it is clear that
	\[
	\frac{s_2}{s_1}=\frac{t_2}{t_1}=1,
	\]
	which means that $s_1=s_2$ and $t_1=t_2$.\\

Finally, we prove that the function $\mathcal{J}$ reaches its unique maximum at the pair $(s_u,t_u)$.
	We know from the above that $(s_u,t_u)$ is the unique critical point of
	$\mathcal{J}$ on $\mathbb{R}_+\times\mathbb{R}_+$. By  $(A_2)$ and  $(A_4)$ , we can find positive constants $C_5$ and $C_6$  such that
	\begin{align}\label{eq:3.18}
	F(x,t) \geq C_{5} \vert t\vert^\theta - C_{6}|t|^{N}, \quad \forall (x,t)  \in (\mathbb{R}^{N},\mathbb{R}\setminus \{0\}).
	\end{align} Using \eqref{eq:1.13} and \eqref{eq:3.18}, we obtain
	\begin{align*}
	\mathcal{J}(s,t)
	&=\mathcal{E}(s u^++t u^-)\\
	&=\frac{1}{N} M(\Vert s u^++t u^- \Vert^N) - \int_{\mathbb{R}^{N}}F(x,s u^+ + t u^-)~dx\\
	&\leq \frac{m(1)}{N} \Vert s u^++t u^- \Vert^N+  \frac{m(1)}{2N} \Vert s u^++t u^- \Vert^{2N} - C_5 |s|^\theta \int_{\mathbb{R}^{N}} |u^+|^\theta dx + C_6 |t|^\theta \int_{\mathbb{R}^{N}} |u^-|^\theta dx\\ &- C_{5}|s|^{\theta}\int_{\mathbb{R}^{N}} |u^+|^\theta dx+C_6 |t|^N\int_{\mathbb{R}^{N}} |u^-|^N dx \\	&\leq \frac{m(1)}{N} (s^N\Vert u^+\Vert^N + t^N\Vert u^-\Vert^N ) + \frac{m(1)}{2N} (s^{2N}\Vert u^+\Vert^{2N} + t^{2N}\Vert u^-\Vert^{2N} )\\ &- C_5 |s|^\theta \int_{\mathbb{R}^{N}} |u^+|^\theta dx - C_5 |t|^\theta \int_{\mathbb{R}^{N}} |u^-|^\theta dx  dx|\\&+C_6 |t|^N\int_{\mathbb{R}^{N}} |u^+|^N dx+C_6 |t|^N\int_{\mathbb{R}^{N}} |u^-|^N dx,
	\end{align*}
which implies that $\displaystyle\lim_{|(s,t)|\to\infty}\mathcal{J}(s,t)=-\infty$, because $\theta>2N$. Hence, it suffices  to show that the maximum point cannot be achieved on the boundary of $\mathbb{R}_+\times\mathbb{R}_+$. We argue by contradiction. Assuming $(0,\bar t)$ is the global maximum point of $\mathcal{J}$ with $\bar t\ge0$,  we have
	\begin{align*}
	\mathcal{J}(s,\bar t)
	&= \frac{1}{N} M(\Vert s u^+ +\bar{t} u^- \Vert^N) - \int_{\mathbb{R}^{N}}F(x,s u^+ + \bar{t} u^-)~dx.
	\end{align*}
	Hence, by \eqref{eq:2.8} it is clear that
		\begin{align*}
	\mathcal{J}_s'(s,\bar t)
	&=m(\|s u^{+} +\bar{t}u^{-}\|^{N})s^{N-1}\|u^{+}\|^{N} - \int_{\mathbb{R}^{N}}f(x,s u^+)u^+~dx >0,
	\end{align*}
	for small enough $s$. This means
	that  $\mathcal{J}$ is an increasing function with
	respect to $s$ if $s$ is small enough, which is a contradiction.
	In a similar way, we can deduce that $\mathcal{J}$ cannot achieve its global maximum at
	$(s,0)$ with $s\ge0$. Thus, we have completed the proof.
\end{proof}

\begin{lemma}\label{lem2.2}
Suppose that $m$ satisfies  $(M_1)-(M_2)$ and assume that $f$ verifies $(A_1)-(A_4)$. Then for any $u \in E$ with $u^{\pm} \neq 0$ such that
	$\langle \mathcal{E}'(u),u^\pm\rangle\leq0$, the unique maximum point pair of $\mathcal{J}$ on $\mathbb{R}_+\times\mathbb{R}_+$ verifies $0<s_u, t_u\leq 1$.
\end{lemma}

\begin{proof}
	Without loss of generality, we may suppose that $ 0< t_u\leq s_u$.
	Since $s_uu^++t_uu^-\in\mathcal{N}$, from $(M_{1})$, we have
	\begin{equation} \label{eq:2.18}
	\begin{aligned}
m(s_u^N \Vert u\Vert^N)s_u^N \Vert u^+\Vert^N &\geq m(s_u^N \Vert u^{+}\Vert^N+ t_u^N \Vert u^{-}\Vert^N)s_u^N \Vert u^+\Vert^N \\&= \int_{\mathbb{R}^{N}}f(x,s_u u^+) s_u u^+ dx.
	\end{aligned}
	\end{equation}
	Furthermore, since $\langle \mathcal{E}'(u),u^+\rangle\leq0$, we have
	\begin{equation}\label{eq 1.19} m(\Vert u\Vert^N) \Vert u^+\Vert^N \leq  \int_{\mathbb{R}^{N}}f(x, u^+) u^+ dx.\end{equation}

Then, from  \eqref{eq:2.18}  and \eqref{eq 1.19}, we get
	
	\begin{equation} \label{eq:2.19}
	\begin{aligned}
	&\frac{m(s_u^N\Vert u\Vert^N)  }{s_u^N \Vert u \Vert^N}
-\frac{m(\Vert u\Vert^N)}{\Vert u\Vert^N} \geq \frac{1}{ \|u\|^{N}\|u^{+}\|^{N}}\int_{\mathbb{R}^{N}}\Big(\frac{f(x,s_u u^+)}{(s_{u} u^+)^{2N-1}}-\frac{f(x,u^+)}{(u^+)^{2N-1}}\Big)(u^+)^{2N}\,dx.
	\end{aligned}
	\end{equation}
From ($A_3$) and ($M_2$), the left hand side of \eqref{eq:2.19} is negative for
$s_u>1$ whereas the right hand side is positive, which is a contradiction. Therefore $0<s_u,t_u\leq 1$.
\end{proof}
\begin{lemma}\label{lem8} Suppose that the hypotheses ($A_{1}$), ($A_{2}$) and ($A_{3}$) are satisfied. Then, for each $x\in B$, we have
	$$tf(x, t) -2N F(x, t)~~ \mbox{is increasing for} ~ t> 0 ~ \mbox{ and decreasing for} ~ t< 0.  $$
	In particular, $tf(x, t) - 2N F(x, t)> 0~~\mbox{for all} ~~(x,t) \in B \times \R\setminus \{0\}.$
\end{lemma}
\begin{proof}
To prove this lemma, we just need to analyse the derivative of $tf(x, t) - 2NF(x, t)$ together with the hypotheses $(A_1)$ and $(A_3)$.
	
\end{proof}
 In the next result, we prove that  sequences in $\mathcal{N}$ cannot converge to $0$ and that $\mathcal{E}$ restricted to $\mathcal{N}$ is bounded from below.
\begin{lemma}\label{lem2.3}
Assume that $m$ satisfies $(M_1)$ and $(M_2)$ and $f$ verifies $(A_2)$ and $(A_4)$. Then, we have
\begin{description}
	\item[$i)$] There exists $\kappa > 0$ such that $\Vert u^+\Vert,$ $\Vert u^-\Vert \geq \kappa,$\: \mbox{ for all}\: $u \in \mathcal{N}$.
	\item[$ii)$] $\displaystyle \mathcal{E}(u)\geq (\frac{1}{2N} - \frac{1}{\theta})m_{0}\|u\|^{N},$ \: \mbox{ for all}\: $u \in \mathcal{N}$.
	\end{description}
\end{lemma}

\begin{proof}\begin{description}
	\item[$i)$] We only prove that there exists $\kappa > 0$ such that $\Vert u^+\Vert \geq \kappa$ for all $u \in \mathcal{N}$ and the proof for $\Vert u^-\Vert$ is similar.
	By contradiction, we suppose that there exists a sequence $\{u_n^+\}\subset \mathcal{N} $ such that $ \Vert u_n^+\Vert \rightarrow 0$ as $ n \rightarrow \infty$.
	Since $u_n \in\mathcal{N}$, we have $\langle I'(u_n),u_n^+\rangle=0$. Thus, by $(M_{1})$ and \eqref{eq:2.7}, we get
	\begin{equation} \label{eq:2.20}
	\begin{aligned}
	 m_{0}\|u^{+}_{n}\|^{N}&\leq m(\Vert u_{n}\Vert^N )\Vert u_{n}^+\Vert^N = \int_{\mathbb{R}^{N}}f(x,u_n^+)~ u_n^+~dx\\ &\leq \e
	\int_{\mathbb{R}^{N}}| u_n^+|^N\,dx + C \int_{\mathbb{R}^{N}}|u_n^+|^{q} exp(\alpha (u_n^+)^{\gamma}-1)dx,
\end{aligned}
	\end{equation}
	for all $n \in \mathbb{N},$ $q>N$ and $\alpha >\alpha_0$. Since $ \Vert u_n^+\Vert \rightarrow 0$ as $ n \rightarrow \infty$, then there exists $n_0 \in \mathbb{N}$ such that $\Vert u_n^+\Vert^\gamma \leq  (\frac{\alpha_{N,\beta}}{2\alpha})$ for all $n \geq n_0$. From H\^older inequality and \eqref{eq:1.7}, we get
	\begin{equation} \label{eq:2.21}
	\begin{aligned}\nonumber
\int_{\mathbb{R}^{N}}|u_n^+|^{q} exp(\alpha (u_n^+)^{\gamma}-1)dx &\leq
 \big(\int_{\mathbb{R}^{N}}|u_n^+|^{2q}dx\big)^{\frac{1}{2}}\big(\int_{\mathbb{R}^{N}} exp(2 \alpha  \|u_{n}^+\|^{\gamma} (\frac{u_{n}^+ }{\Vert u_{n}^+\Vert})^{\gamma}-1)dx\big)^{\frac{1}{2}}
	\\&\leq  C \Vert u_n^+\Vert_{2q}^q.
	\end{aligned}
	\end{equation}
	Combining \eqref{eq:2.20} with the last inequality, we can deduce from the Sobolev embedding Theorem that when $n \geq n_0$,
\begin{equation} \label{eq:2.22}
\begin{aligned}
m_{0}\|u_n^+\|^{N} \leq C_{3} \e \Vert u_n^+\Vert^N + C_{4} \Vert u_n^+\Vert^q.
\end{aligned}
\end{equation}
We can choose $\e > 0$ such that $(m_{0}- C_{3} \e) > 0$ and since $q>N$, we can deduce that \eqref{eq:2.22} contradicts  $ \Vert u_n^+\Vert \rightarrow 0$ as $ n \rightarrow \infty$. The proof is complete.

\item[$ii)$] Given $u \in \mathcal{N}$, by the definition of $\mathcal{N}$, $(M_{1})$, (\ref{eq:1.15}) and $(A_2)$, we obtain
	\begin{equation} \label{eq:2.23}
	\begin{aligned}
	\mathcal{E}(u)
	& = \mathcal{E}(u)-\frac{1}{\theta}\langle \mathcal{E}'(u),u\rangle \\
	& =\displaystyle\frac{1}{N}M(\|u\|^{N}) - \frac{1}{\theta}m(\| u\|^{N})\| u\|^{N}  - \frac{1}{\theta}\int_{\mathbb{R}^{N}}\big(f(x, u) u- \theta F(x, u)\big)dx \\
	& \ge(\frac{1}{2N} - \frac{1}{\theta})m_{0}\|u\|^{N}.
	\end{aligned}
	\end{equation}
		\end{description}
	
\end{proof}

	So, we have $\mathcal{E}(u)>0$ for all $u\in\mathcal{N}$. Therefore, $\mathcal{E}(u)$ is bounded
	below on $\mathcal{N},$ that is  $\displaystyle c_{\mathcal{N}}=\inf_{u\in\mathcal{N}}\mathcal{E}(u)$ is well-defined.\\
To end this section, let us prove that if the minimum of $\mathcal{E}$ over $\mathcal{N}$  is reached at some $u \in \mathcal{N}$, then $u$ is a critical point of $\mathcal{E}$.
\begin{lemma}\label{lemma4}
If $u_0\in \mathcal{N}$ satisfies $\mathcal{E}(u_0)= c_{\mathcal{N}} $, then $ \mathcal{E}'(u_0)= 0.$
	\end{lemma}
\begin{proof}
Suppose by contradiction that $\mathcal{E}'(u_0)\neq0$. By the continuity of $\mathcal{E}'$, it follows that there exist $\gamma>0$ and
	$\iota>0$ such that
	\[
	\|\mathcal{E}'(v)\|\geq\iota \quad \text{for all}\,\,\|v-u_0\|
	\leq 3\gamma.
	\]
	Choose $\nu \in(0,\min\{1/2,\frac{\gamma}{\sqrt{2}\|u_0\|}\})$. Let
	$
	D:= (1-\nu,1+\nu)\times(1-\nu,1+\nu)
	$ and $$
	k(s,t):= s u_0^+ + t u_0^-\quad \text{for all}\, \, (s,t)\in D.
	$$
	In view of Lemma \ref{nonempty}, we have
	\begin{equation}\label{eq:2.24}
	\displaystyle\bar{c}_{\mathcal{N}} :=\max_{\partial D}(\mathcal{E} \circ k)< c_{\mathcal{N}}.
	\end{equation}
	Let $\varepsilon:= \min\{(c_{\mathcal{N}}
	-\bar{c}_{\mathcal{N}})/3,\iota\gamma/8\}$ and Denote by $S_{r}:=B(u_{0},r),r>0$ the ball centered at $ u_{0}$ and of radius $r$.
	According to Lemma 2.3 in \cite{mw1996},
	there exists a deformation $\eta\in C([0,1]\times E,E)$  such that
	\begin{itemize}
		\item [(a)] $\eta(1,v) = v$ if $v\notin(\mathcal{E}^{-1}([c_{\mathcal{N}} -2\varepsilon, c_{\mathcal{N}} +2\varepsilon])\cap
		S_{2\gamma})$,
		
		\item[(b)] $\eta(1,\mathcal{E}^{c_{\mathcal{N}} +\varepsilon}\cap
		S_{\gamma}) \subset \mathcal{E}^{c_{\mathcal{N}} -\varepsilon}$,
		
		\item[(c)] $\mathcal{E}(\eta(1,v))\leq \mathcal{E}(v)$ for all $v\in E$.
	\end{itemize}
	
	Clearly,
	\begin{equation}\label{eq:2.25}
	\max_{(s,t)\in\bar{D}}\mathcal{E}(\eta(1,k(s,t)))< c_{\mathcal{N}}.
	\end{equation}
	Therefore we claim that $\eta(1,k(D))\cap
	\mathcal{N} \neq \emptyset$ , which  contradicts the definition of $c_{\mathcal{N}}$.
	
	We define $\bar{k}(s,t):=\eta(1,k(s,t))$,
	\begin{align*}
	\Gamma_0(s,t)
	:&=(\langle \mathcal{E}'(k(s,t)), su_0^+\rangle,\langle \mathcal{E}'(k(s,t)), t u_0^-\rangle)\\
	& = (\langle \mathcal{E}'(s u_0^+ + t u_0^-),
	s u_0^+\rangle, \langle \mathcal{E}'(s u_0^+ + t u_0^-), t u_0^-\rangle)\\
	& = ( \phi^1_u(s,t) , \phi^2_u(s,t) )
	\end{align*}
	and
	\[
	\Gamma_1(s, t) := \Big(\frac{1}{s}\langle \mathcal{E}'(\bar{k}(s,t)),\,
	(\bar{k}(s,t))^+\rangle,\frac{1}{t}\langle \mathcal{E}'(\bar{k}(s,t)),(\bar{k}(s,t))^-\rangle\Big).
	\]
	
By a straightforward computation and using (\ref{eq:1.14}), we get	
	
\begin{align}\label{eq:2.26}
\frac{\del \phi_u^1(s,t)}{\del s}|_{(1,1)}\nonumber
:&= Nm'(\|u_{0}\|^{N})\|u^{+}_{0}\|^{2N}+Nm(\|u_{0}\|^{N})\|u^{+}_{0}\|^{N}\\&\nonumber - \int_{\mathbb{R}^{N}}(f'(x,u^{+}_{0}) (u^{+}_{0})^2 + f(x,u^{+}_{0}) u^{+}_{0}) dx\\ &<2Nm(\|u_{0}\|^{N})\|u^{+}_{0}\|^{N}-Nm'(\|u_{0}\|^{N})\|u^{+}_{0}\|^{N}\|u^{-}_{0}\|^{N}\\&\nonumber - \int_{\mathbb{R}^{N}}(f'(x,u^{+}_{0}) (u^{+}_{0})^2 + f(x,u^{+}_{0}) u^{+}_{0}) dx.
\end{align}	
\begin{align}\nonumber
\frac{\del \phi_u^2(s,t)}{\del t}|_{(1,1)}
&= Nm'(\|u_{0}\|^{N})\|u^{-}_{0}\|^{2N}+Nm(\|u_{0}\|^{N})\|u^{-}_{0}\|^{N}\\&-\nonumber \int_{\mathbb{R}^{N}}(f'(x,u^{-}_{0}) (u^{-}_{0})^2 + f(x,u^{-}_{0}) u^{-}_{0}) dx\\\nonumber
& <2Nm(\|u_{0}\|^{N})\|u^{-}_{0}\|^{N}-Nm'(\|u_{0}\|^{N})\|u^{+}_{0}\|^{N}\|u^{-}_{0}\|^{N}\\&\nonumber-
\int_{\mathbb{R}^{N}}(f'(x,u^{-}_{0}) (u^{-}_{0})^2 + f(x,u^{-}_{0}) u^{-}_{0}) dx
\end{align}
and
\begin{equation*}
\frac{\del \phi_u^1(s,t)}{\del t}|_{(1,1)} = \frac{\del \phi_u^2(s,t)}{\del s}|_{(1,1)}
:= Nm'(\|u_{0}\|^{N})\|u^{+}_{0}\|^{N}\|u^{-}_{0}\|^{N}.
\end{equation*}	
On the other hand, since $\langle \mathcal{E}'(u_{0}),u_{0}^+\rangle=0$, we have
\begin{equation}\label{eq:2.27}
m(\|u_{0}\|^{N})\|u^{+}_{0}\|^{N}= \int_{\mathbb{R}^{N}}f(x,u^{+}_{0})~ u^{+}_{0}~dx.
\end{equation}
Combining \eqref{eq:2.26} and \eqref{eq:2.27} with the condition $(H_3)$ and $u^{+}\ne 0$, we get
\begin{align}\label{eq:2.28}
\frac{\del \phi_u^1(s,t)}{\del s}|_{(1,1)}& \nonumber
<-Nm'(\|u_{0}\|^{N})\|u^{+}_{0}\|^{N}\|u^{-}_{0}\|^{N} - \nonumber\int_{\mathbb{R}^{N}}(f'(x,u^+) (u^+)^2 - (2N-1) f(x,u^+) u^+) dx\\ &<-Nm'(\|u_{0}\|^{N})\|u^{+}_{0}\|^{N}\|u^{-}_{0}\|^{N}=-\frac{\del \phi_u^1(s,t)}{\del t}|_{(1,1)}<0 .
\end{align}

Similarly, we can show that
\begin{equation}\label{eq:2.29}
\frac{\del \phi_u^2(s,t)}{\del t}|_{(1,1)}  <- \frac{\del \phi_u^1(s,t)}{\del t}|_{(1,1)}  <0.
\end{equation}

Let
$$\displaystyle K=\begin{bmatrix} \frac{\del \phi_u^1(s,t)}{\del s}|_{(1,1)} & \frac{\del \phi_u^2(s,t)}{\del s}|_{(1,1)} \\ \\ \frac{\del \phi_u^1(s,t)}{\del t}|_{(1,1)} & \frac{\del \phi_u^2(s,t)}{\del t}|_{(1,1)} \end{bmatrix}. $$	
Hence, by \eqref{eq:2.28} and \eqref{eq:2.29}, we have
\begin{align*}\det K &= \frac{\del \phi_u^1(s,t)}{\del s}|_{(1,1)} \times \frac{\del \phi_u^2(s,t)}{\del t}|_{(1,1)} - \frac{\del \phi_u^2(s,t)}{\del s}|_{(1,1)} \frac{\del \phi_u^1(s,t)}{\del t}|_{(1,1)}\\
&= \frac{\del \phi_u^1(s,t)}{\del s}|_{(1,1)} \times \frac{\del \phi_u^2(s,t)}{\del t}|_{(1,1)} -  \Big(\frac{\del \phi_u^1(s,t)}{\del t}|_{(1,1)}\Big)^2 \\
&> \Big(\frac{\del \phi_u^1(s,t)}{\del t}|_{(1,1)}\Big)^2 -  \Big(\frac{\del \phi_u^1(s,t)}{\del t}|_{(1,1)}\Big)^2 = 0.
\end{align*}

Therefore, $\Gamma_0(s,t)$ is a $C^1$ function and $(1, 1)$ is the unique isolated zero point of $\Gamma_0$. By using the degree theory, we deduce that 	$\deg(\Gamma_0 ,D,0) = 1$.
	
Hence, combining \eqref{eq:2.24} with $(a)$, we obtain
	\[
	k(s,t)
	= \bar{k}(s,t)\quad \text{on } \partial D.
	\]
	Therefore, by the degree theory (see [\cite{DV}, Theorem 4.5]), we get $\deg(\Gamma_1 ,D,0) =\deg(\Gamma_0 ,D,0) = 1$.
	Hence, again by the degree theory $\Gamma_1(s_0, t_0) = 0$ for some $(s_0, t_0)\in D$
	so that
	\[
	\eta(1,k(s_0, t_0))=\bar{k}(s_0, t_0)\in \mathcal{N},
	\]
	which contradicts \eqref{eq:2.25}. Hence, $\mathcal{E}'(u_0)=0$,
	which implies that $u_0$ is a critical point of $\mathcal{E}$.
\end{proof}

\section{Proof of Theorem \ref{th1}}
\begin{lemma} \label{lemma7} There exists $w \in \mathcal{N}$ such that $\mathcal{E}(w) = c_{\mathcal{N}}$.
\end{lemma}

\begin{proof}
	Let sequence $(w_n) \subset \mathcal{N} $ satisfy $\displaystyle \lim_{n \rightarrow +\infty} \mathcal{E}(w_n) = c_{\mathcal{N}}$. It is clearly that $(w_n)$ is bounded by Lemma \ref{lem2.3}. Then, up to a subsequence, there exists $w \in E$ such that
	
	\begin{equation}\label{eq:3.8}\begin{array}{ll}
	w_{n}^\pm \rightharpoonup w^\pm~~~~&\mbox{in}~~E,\\
	w_{n}^\pm \rightarrow w^\pm~~&\mbox{in}~~L^{q}(\mathbb{R}^{N}),~~\forall q\geq N,\\
	w_{n}^\pm \rightarrow w^\pm ~~&\mbox{a.e. in }~~\mathbb{R}^{N}.
	\end{array}
	\end{equation}
	We claim that
	\begin{equation}\label{eq:3.7}
	\int_{\mathbb{R}^{N}} f(x, w_n^\pm)~ w_n^\pm dx \rightarrow \int_{\mathbb{R}^{N}} f(x, w^\pm)~ w^\pm dx.
	\end{equation}
	Indeed, by \eqref{eq:2.7}, we have
	\begin{equation}\label{eq:3.2}
	\int_{\mathbb{R}^{N}}f(x, w_n^\pm)~ w_n^\pm dx \leq \e
	\int_{\mathbb{R}^{N}}| w_n^\pm|^N\,dx + C \int_{\mathbb{R}^{N}}|w_n^\pm|^{q} exp(\alpha |w_n^\pm|^{\gamma})dx, \quad \mbox{for all}\; \alpha > 0\, \mbox{and}\, q>N.
	\end{equation}
	We define $g(w_n^\pm(x))$ as follows
	\begin{equation}\label{eq3.2} g(w_n^\pm(x)) := \e|w_n^\pm|^N + C|w_n^\pm|^{q} exp(\alpha |w_n^\pm|^{\gamma}).\end{equation}
	We will prove that $g(w_n^\pm(x))$ is convergent in $L^1(B)$. First note that
		\begin{equation}\label{eq:3.3'}
	|w_{n}|^N \rightarrow |w|^N~~ \mbox{in}~~ L^{1}(\mathbb{R}^{N}).
	\end{equation}
	Considering $s,s' > 1$ such that $\frac{1}{s} +\frac{1}{s'}=1,$ we get
	\begin{equation}\label{eq:3.3}
	|w_{n}|^q \rightarrow |w|^q~~ \mbox{in}~~ L^{s'}(\mathbb{R}^{N}).
	\end{equation}
	Moreover, choosing $\alpha>0$ enough small such that $\displaystyle s\alpha (\max_{n} ||w_n^\pm||^{\gamma}) \leq \alpha_{N,\beta},$ we conclude from Theorem \ref{th3} that
	\begin{equation}\label{eq:3.4}
	\int_{\mathbb{R}^{N}} exp(s \alpha |w_n^\pm|^{\gamma})dx \leq M.
	\end{equation}
	Since
	\begin{equation}\label{eq:3.5}
	exp(\alpha |w_n^\pm|^{\gamma}) \rightarrow exp(\alpha |w^\pm|^{\gamma})\:  \mbox{a.e in }~~B.
	\end{equation}
	Then, from \eqref{eq:3.4} and [\cite{Ka}, Lemma 4.8, Chapter 1], we get that
	\begin{equation}\label{eq:3.6}
	exp(\alpha |w_n^\pm|^{\gamma}) \rightharpoonup exp(\alpha |w^\pm|^{\gamma})\:  \mbox{ in }~~L^{s}(B).
	\end{equation}
	
	Then it follows from the H\"older inequality, \eqref{eq:3.3'}, \eqref{eq:3.3}, \eqref{eq:3.6} and Trudinger-Moser inequality that
	\begin{align}\nonumber
	&\int_{\mathbb{R}^{N}}\big(g(w_n^\pm(x)) - g(w^\pm(x))\big)dx= \e \int_{\mathbb{R}^{N}}\big(|w_n^\pm|^{N} - |w^\pm|^{N}\big) dx \\& + C \int_{\mathbb{R}^{N}}\big(|w_n^\pm|^{q} - |w^\pm|^{q}\big) exp(\alpha |w_n^\pm|^{\gamma})dx + C \int_{\mathbb{R}^{N}}|w^\pm|^{q} \big(exp(\alpha |w_n^\pm|^{\gamma}) - exp(\alpha |w^\pm|^{\gamma})\big) dx\nonumber \\&\leq\e \int_{\mathbb{R}^{N}}\big(|w_n^\pm|^{N} - |w^\pm|^{N}\big) dx+ \nonumber C \Big(\int_{\mathbb{R}^{N}}\big(|w_n^\pm|^{q} - |w^\pm|^{q}\big)^{s'}dx \Big)^{\frac{1}{s'}} \Big( \int_{\mathbb{R}^{N}} exp(s\alpha |w^\pm|^{\gamma})dx\Big)^{\frac{1}{s}}\\&\quad + \nonumber C \int_{\mathbb{R}^{N}}|w^\pm|^{q} \big(exp(\alpha |w_n^\pm|^{\gamma}) - exp(\alpha |w^\pm|^{\gamma})\big) dx\\
	& \leq \nonumber \e \int_{\mathbb{R}^{N}}\big(|w_n^\pm|^{N} - |w^\pm|^{N}\big) dx+CM\Big(\int_{\mathbb{R}^{N}}\big(|w_n^\pm|^{q} - |w^\pm|^{q}\big)^{s'}dx \Big)^{\frac{1}{s'}}  \\ \nonumber&+ C \int_{\mathbb{R}^{N}}|w^\pm|^{q} \big(exp(\alpha |w_n^\pm|^{\gamma}) - exp(\alpha |w^\pm|^{\gamma})\big) dx\\& \nonumber \rightarrow 0 \: \mbox{as} \: n \rightarrow \infty,
	\end{align}
	which closes the proof of \eqref{eq:3.7}. Thus, as a direct consequence,
	\begin{equation}\label{eq 3.8}\int_{\mathbb{R}^{N}}F(x,w_{n}^\pm) dx \rightarrow \int_{\mathbb{R}^{N}}F(x,w^\pm)dx, \end{equation}
	hold.\\
Now, we claim that $w^\pm\ne 0$. Suppose, by contradiction, $w^+= 0$. From the definition
	of $\mathcal{N}$, \eqref{eq:3.8}  and \eqref{eq:3.7}, we have that $ \displaystyle\lim_{n \rightarrow +\infty} \Vert w^+_n \Vert =0$, which contradicts Lemma \ref{lem2.3}. Hence, $w^+\ne 0$ and $w^-\ne 0$.
	
	From the lower semicontinuity of norm and \eqref{eq:3.8}, it follows that
 	\begin{equation}\label{eq:3.9}
 	m(\Vert w\Vert^N)\Vert w^+\Vert^N \leq \liminf_{n \rightarrow +\infty}m(\Vert w_{n}\Vert^N)\Vert w_{n}^+\Vert^N.
 	\end{equation}	

 On the other hand, by using $\langle \mathcal{E}'(w_n),w_n^+\rangle =0$ and \eqref{eq:3.7}, we have

 	\begin{equation}\label{eq:3.10} \liminf_{n \rightarrow +\infty}m(\Vert w_{n}\Vert^N)\Vert w_{n}^+\Vert^N= \liminf_{n \rightarrow +\infty} \int_{\mathbb{R}^{N}}f(x, w_n^+) w_n^+ dx=\int_{\mathbb{R}^{N}}f(x, w^+) w^+ dx.\end{equation}
 From \eqref{eq:3.9} and \eqref{eq:3.10} we deduce that $\langle \mathcal{E}'(w),w^+\rangle \leq 0$, and similarly we can prove $\langle \mathcal{E}'(w),w^-\rangle \leq 0$.
 Then, Lemma \ref{lem2.2} implies that there exists $(s_u, t_u) \in (0, 1] \times (0, 1] $ such that $s_uw^++t_uw^-\in\mathcal{N}$. Thus, by the lower semicontinuity of norm, (\ref{eq:1.15}), \eqref{eq:3.7}, \eqref{eq 3.8} and Lemma \ref{lem8}, we get that

$$\begin{array}{rclll}
\displaystyle c_{\mathcal{N}} \leq \mathcal{E}(s_uw^++t_uw^-)&=&\mathcal{E}(s_uw^++t_uw^-) -\displaystyle\frac{1}{2N} \langle \mathcal{E}'(s_uw^++t_uw^-),s_uw^++t_uw^-\rangle
\\ &\leq& \mathcal{E}(w) -\displaystyle\frac{1}{2N} \langle \mathcal{E}'(w),w\rangle\\
&=&\displaystyle \frac{1}{N}M(\|w\|^{N})-\frac{1}{2N}m(\|w\|^{N})\|w\|^{N}~\\ &+&\displaystyle\frac{1}{2N}\int_{\mathbb{R}^{N}}\big(f(x,w) w-2NF(x,w)\big)dx\\
&\leq&\displaystyle\liminf_{n\rightarrow +\infty}\Big[\displaystyle\frac{1}{N}M(\|w_{n}\|^{N})-\frac{1}{2N}m(\|w_{n}\|^{N})\|w_{n}\|^{N}\\&+&\displaystyle\frac{1}{2N}\int_{\mathbb{R}^{N}}\big(f(x, w_{n}) w_n-2NF(x, w_{n})\big)dx\Big]\\
&\leq&\displaystyle \liminf_{n\rightarrow+\infty}\big[\mathcal{E}(w_n) -\displaystyle\frac{1}{2N} \langle \mathcal{E}'(w_n), w_n\rangle\big] = c_{\mathcal{N}}.
\end{array}$$	
Therefore, we get that $\mathcal{E}(w)=c_{\mathcal{N}}$, which is the desired conclusion.	
\end{proof}
Thus, from Lemma \ref{lemma4},  $w$ is a least energy sign-changing solution of problem \eqref{eq:1.1}.

Now, in the case when $m(t)=a+bt, a, b>0$, we prove that $w$ has exactly two nodal domains. To this end, we prove the following Lemma.
\begin{lemma}\label{lemma8}Suppose that $m(t)=a+bt, a>0 ,b>0$. Then,
	if $w$ is a least energy sign-changing solution of problem \eqref{eq:1.1}, then $w$ has
	exactly two nodal domains.
\end{lemma}
\begin{proof} Assume by contradiction
that $w= w_1 + w_2 + w_3$ satisfies
$$ w_i \ne 0, i=1,2,3, w_1 \geq 0, w_2 \leq 0, \;\mbox{a.e. in}\; \mathbb{R}^{N} $$
$$\Omega_1 \cap \Omega_2= \emptyset,  \Omega_1 = \{ x \in \mathbb{R}^{N} :w_1(x) > 0 \},\; \Omega_2 = \{ x \in \mathbb{R}^{N}  : w_2 (x) < 0 \} $$
$$w_1|_{\mathbb{R}^{N} \setminus \Omega_1 \cup \Omega_2} = w_2|_{ \mathbb{R}^{N}\setminus \Omega_2 \cup \Omega_1}= w_3|_{ \Omega_1 \cup \Omega_2} =0 $$
and
\begin{equation}\label{eq:3.11}
\langle \mathcal{E}'(w), w_i\rangle =0 \:\mbox{for}\: i=1,2,3.
\end{equation}
Setting $v = w_1 + w_2,$ we have that $v^+ = w_1$ and $v^- = w_2,$ i.e. $v^\pm \ne 0$. From Lemma \ref{nonempty}, it follows that there exists a unique point pair $(s_v, t_v) \in \mathbb{R}_+\times \mathbb{R}_+$  such that
$s_v w_1 + t_v w_2 \in \mathcal{N}.$ Hence,
$\mathcal{E}(s_v w_1 + t_v w_2) \geq c_{\mathcal{N}}.$
Moreover, using the fact that $\langle \mathcal{E}'(w), w_i\rangle =0$, we obtain $$\langle \mathcal{E}'(v), v^\pm \rangle = -b \| v^\pm\|^{N}\| w_3\|^{N}  < 0.$$ From Lemma \ref{lem2.2}, we have that $$(s_v, t_v) \in (0, 1] \times (0, 1].$$

On the other hand, by $(A_2)$ we have that
$$\begin{array}{rclll}
\displaystyle 0= \frac{1}{2N} \langle \mathcal{E}'(w), w_3\rangle &=&  \displaystyle\frac{a}{2N}\|w_3\|^{N} + \frac{b}{2N}\| w_1\|^{N}\| w_3\|^{N} + \frac{b}{2N} \| w_2\|^{N}\| w_3\|^{N}  + \frac{b}{2N}\| w_3\|^{2N}\\ && - \displaystyle\frac{1}{2N}\int_{\mathbb{R}^{N}} f(x, w_3)w_3 dx\\
& <&\displaystyle \mathcal{E}(w_3) +  \frac{b}{2N}\| w_1\|^{N}\| w_3\|^{N} + \frac{b}{2N} \| w_2\|^{N}\| w_3\|^{N}.
\end{array}$$

Hence, by Lemma \ref{lem8}, we can obtain that
$$\begin{array}{rclll}
\displaystyle m \leq \mathcal{E}(s_v w_1 +t_v w_2)&=&\mathcal{E}(s_v w_1 +t_v w_2) -\displaystyle\frac{1}{2N} \langle \mathcal{E}'(s_v w_1 +t_v w_2),s_v w_1 +t_v w_2\rangle
\\ &=& \displaystyle a\frac{s_v^N}{2N}\|w_1\|^{N} + a\frac{t_v^N}{2N}\|w_2\|^{N}  +\frac{1}{2N}\int_{\mathbb{R}^{N}}\big(f(x,s_v w_1) s_v w_1-2NF(x,s_v w_1)\big)dx \\&& \displaystyle +\frac{1}{2N}\int_{\mathbb{R}^{N}}\big(f(x,t_v w_2) t_v w_2-2NF(x,t_v w_2)\big)dx
\\ &\leq & \displaystyle\frac{a}{2N}\|w_1\|^{N} + \frac{a}{2N}\|w_2\|^{N} +\frac{1}{2N}\int_{\mathbb{R}^{N}}\big(f(x, w_1) w_1-2NF(x, w_1)\big)dx\\&&\displaystyle +\frac{1}{2N}\int_{\mathbb{R}^{N}}\big(f(x, w_2)  w_2-2NF(x, w_2)\big)dx
\\ &=& \mathcal{E}(w_1 + w_2) -\displaystyle\frac{1}{2N} \langle \mathcal{E}'(w_1 + w_2),w_1 + w_2\rangle\\
&=& \displaystyle \mathcal{E}(w_1 + w_2) + \frac{1}{2N} \langle \mathcal{E}'(w), w_3\rangle + \frac{b}{2N}\| w_1\|^{N}\| w_3\|^{N} + \frac{b}{2N} \| w_2\|^{N}\| w_3\|^{N}
\\
&<&\displaystyle \mathcal{E}(w_1) + \mathcal{E}(w_2) + \mathcal{E}(w_3) + \frac{b}{2N} \| w_1\|^{N}\| w_2\|^{N}  + \frac{b}{2N} \| w_1\|^{N}\| w_3\|^{N}\\&&\displaystyle + \frac{b}{2N} \| w_2\|^{N}\| w_3\|^{N}\\
&=& \mathcal{E}(w) = m,
\end{array}$$	
which is a contradiction, that is, $w_3 =0$ and $w$ has exactly two nodal domains.
\end{proof}
\section{Proof of Theorem \ref{th2}}
\subsection{Auxiliary problem}
Before proving theorem \ref{th2}, we need to consider the following auxiliary problem
\begin{equation}\label{eq:4.1}
m(\|u\|^{N})\big[-\textmd{div} (\omega_{\beta}(x)|\nabla u|^{N-2}  \nabla u)\big] = \displaystyle \vert u\vert^{p-2}u~~ \mbox{in} ~~ \mathbb{R}^{N}
\end{equation}
where $p$ is the constant that appear in the assumption $(A_5)$. The energy functional $\mathcal{E}_p$ associated to \eqref{eq:4.1} is given by
$$\mathcal{E}_p(u) :=\frac{1}{N}M(\|u\|^{N}) -\frac{1}{p}\int_{\mathbb{R}^{N}} \vert u\vert^p dx  $$ and  the sign-changing Nehari set is defined by
$$\mathcal{N}_p:= \{ u \in E, u^{\pm} \neq 0\; \mbox{and}\; \langle \mathcal{E}_p'(u),u^+ \rangle=\langle \mathcal{E}_p'(u),u^-\rangle = 0\}.$$
Let $ c_{\mathcal{N}_p} = \inf_{\mathcal{N}_p} \mathcal{E}_p(u)$, we have the following result.
\begin{lemma}\label{lemma 8}
	There exists $w \in \mathcal{N}_p$ such that $\mathcal{E}_p(w) = c_{\mathcal{N}_p}$.
\end{lemma}
\begin{proof}
	This lemma can be proved in four stages:
	\begin{description}
		\item[$ {\bf Stage\; 1.}$]  For any $u \in E$ with $u^{\pm} \neq 0$, similarly to the Lemma \ref{nonempty}, there is a unique maximum point pair $(s_u, t_u) \in \mathbb{R}_+\times \mathbb{R}_+$ of the function $\mathcal{E}_p$ such that $s_u u^+ + t_u u^- \in \mathcal{N}_p.$
		\item[$ {\bf Stage\; 2.}$] If $u \in E$ with $u^{\pm} \neq 0$, such that
		$\langle \mathcal{E}_p'(u),u^\pm\rangle\leq0$. Then, similar to Lemma \ref{lem2.2}, the unique maximum point pair $(s_u, t_u)$ in Stage (1) satisfies $0<s_u, t_u\leq 1$.
		\item[$ {\bf Stage\; 3.}$]By analogy with the Lemma \ref{lem2.3}, for all $u \in \mathcal{N}_p$, there exists $\kappa > 0$ such that $\Vert u^+\Vert,$ $\Vert u^-\Vert \geq \kappa.$
		\item[$ {\bf Stage\; 4.}$] Now, let a sequence $(w_n) \subset \mathcal{N}_p $ satisfying  $\displaystyle \lim_{n \rightarrow +\infty} \mathcal{E}(w_n) = c_{\mathcal{N}_p}$. Similar to Lemma \ref{lemma7}, we can show that, up to a subsequence, 	$w_{n}^\pm \rightharpoonup w^\pm~~~~\mbox{in}~~E.$ From Stage (3), we show that $w^\pm \ne 0$. Using the Stages (1), (2) and again similar to Lemma \ref{lemma7}, we get $w \in \mathcal{N}_p$ such that $\mathcal{E}_p(w) = c_{\mathcal{N}_p}$.
	\end{description}
\end{proof}
\subsection{Estimation of the nodal level $c_{\mathcal{N}}$}
We will now obtain an important estimate of the nodal level $c_{\mathcal{N}}$. This is a powerful tool for obtaining an appropriate bound on the norm of a minimising sequence for $c_{\mathcal{N}}$ in $\mathcal{N}$.

\begin{lemma}\label{lem9}
	Assume that $(M_1), (M_2),$ $(A_1), (A_5)$ and \eqref{eq:1.18} are satisfied. It holds that
		\begin{equation} \label{eq:4.2} c_{\mathcal{N}} \leq m_{0}\frac{\theta-2N}{2N\theta} \bigg(\frac{\alpha_{N,\beta}}{ 2(\alpha_0 + \delta)}\bigg)^{(N-1)(1-\beta)}.\end{equation}
\end{lemma}

\begin{proof}
From Lemma \ref{lemma 8}, there exists $w \in \mathcal{N}_p$ such that $\mathcal{E}_p(w) = c_{\mathcal{N}_p}$ and $\mathcal{E}_p'(w)=0.$ We therefore obtain

\begin{equation}\label{e2.9'}
\frac{1}{N} M(\Vert w\Vert^N) - \frac{1}{p}\int_{\mathbb{R}^{N}} \vert w\vert^p ~dx =  c_{\mathcal{N}_p}\quad \mbox{and}\quad m( \Vert w\Vert^N) \Vert w^\pm\Vert^N= \int_{\mathbb{R}^{N}} \vert w^\pm\vert^p ~dx.
\end{equation}
Together with \eqref{eq:1.15}, this implies that
\begin{equation}\label{eq11'}
(\frac{1}{2N}- \frac{1}{p})|w|_p^p =\frac{1}{2N}m( \Vert w\Vert^N) \Vert w\Vert^N - \frac{1}{p}|w|_p^p\\
\leq \frac{1}{N} M(\Vert w\Vert^N) - \frac{1}{p}|w|_p^p = c_{\mathcal{N}_p}.
\end{equation}

On the other hand, using $(A_5)$ and \eqref{e2.9'}, we get $\langle \mathcal{E}'(w),w^\pm\rangle \leq 0$ which together with Lemma \ref{lem2.2} yielding that there is a unique pair $(s, t) \in (0, 1] \times (0, 1] $ such that $sw^++tw^-\in\mathcal{N}$. Using $(A_5)$, $(M_1)$, (\ref{eq:1.14}), \eqref{e2.9'} and \eqref{eq11'}  we obtain
$$\begin{array}{rclll}
\displaystyle c_{\mathcal{N}} &\leq& \mathcal{E}(sw^++tw^-)
\\  &\leq& \displaystyle\frac{m(1)}{N}s^N\|w^+\|^N+\frac{m(1)}{2}t^N\|w^-\|^N +  \frac{m(1)}{2N}s^{2N}\|w^+\|^{2N} +
\frac{m(1)}{2N}t^{2N}\|w^-\|^{2N}   \\
&&\displaystyle +\frac{m(1)}{N}s^N t^N\|w^+\|^N\|w^-\|^N - \frac{C_p}{p} s^p\vert w^+\vert_p^p -  \frac{C_p}{p} t^p\vert w^-\vert^p_p \\
&\leq&  \displaystyle\frac{m(1)}{Nm_{0}}s^N |w^{+}|_p^p+\frac{m(1)}{Nm_{0}}t^N |w^{-}|_p^p +  \frac{m(1)}{2Nm^{2}_{0}}s^{N} |w^{+}|^{2p}_{p} +
\frac{m(1)}{2Nm^{2}_{0}}t^{N}|w^{-}|^{2p}_{p}  \\
&&\displaystyle +\frac{m(1)}{Nm_0^2}s^N t^N|w^{+}|_p^p|w^{-}|_p^p - \frac{C_p s^p}{p}\vert w^+\vert_p^p -  \frac{C_p t^p}{p}\vert w^-\vert^p_p \\
&\leq&\displaystyle \max_{\xi >0}((\frac{m(1)}{m_{0}}+\frac{m(1)}{2m^{2}_{0}}|w|_p^p)\frac{\xi^N}{N} - C_p\frac{\xi^p}{p}) |w|_p^p \\ &\leq&\displaystyle \max_{\xi >0}\big((\frac{m(1)}{m_{0}}+\frac{m(1)}{m^{2}_{0}}\frac{pN}{p-2N}c_{\mathcal{N}_{p}})\frac{\xi^N}{N} - C_p\frac{\xi^p}{p}\big) \frac{2pN}{p-2N}c_{\mathcal{N}_{p}}.

\end{array}$$
With a few simple algebraic manipulations, we obtain
\begin{equation}\label{eq:4.7}
c_{\mathcal{N}} \leq  (\frac{m(1)}{m_{0}}+\frac{m(1)}{m^{2}_{0}}\frac{pN}{p-2N}c_{\mathcal{N}_{p}})^{\frac{p}{p-N}}C_p^{\frac{-N}{p-N}}\frac{2(p-N)}{p-2N}c_{\mathcal{N}_{p}}.
\end{equation}
Therefore, by \eqref{eq:1.18} and \eqref{eq:4.7}, we obtain that \eqref{eq:4.2} holds.
\end{proof}

The next result gives us some campactness properties of minimizing sequences.
\begin{lemma} \label{lemma7'} \begin{itemize}\item [(i)]If  $(w_n) \subset \mathcal{N} $ is a minimizing sequence for $ c_{\mathcal{N}}$, then up to a subsequence  there exists $w \in E$ such that
			\begin{equation*}\label{eq:3.81}
		w_{n}^\pm \rightharpoonup w^\pm~~~~\mbox{in}~~E,\:
		w_{n}^\pm \rightarrow w^\pm~~\mbox{in}~~L^{q}(\mathbb{R}^{N}),~~\forall q\geq N\; \mbox{and}\:
		w_{n}^\pm \rightarrow w^\pm ~~\mbox{a.e. in }~~\mathbb{R}^{N}.
		\end{equation*}
	and
		\begin{equation}\label{eq:3.71}
		\int_{\mathbb{R}^{N}} f(x, w_n^\pm)~ w_n^\pm dx \rightarrow \int_{\mathbb{R}^{N}} f(x, w^\pm)~ w^\pm dx.
		\end{equation}
\item[(ii)] There exists $w \in \mathcal{N}$ such that $\mathcal{E}(w) = c_{\mathcal{N}}$.
\end{itemize}
\end{lemma}
\begin{proof}$(i)$	
		Let sequence $(w_n) \subset \mathcal{N} $ satisfy $\displaystyle \lim_{n \rightarrow +\infty} \mathcal{E}(w_n) = c_{\mathcal{N}}$. It is clearly that $(w_n)$ is bounded by Lemma \ref{lem2.3}. Then, up to a subsequence, there exists $w \in E$ such that
			\begin{equation*}\label{eq:3.81}
		w_{n}^\pm \rightharpoonup w^\pm~~~~\mbox{in}~~E,\:
		w_{n}^\pm \rightarrow w^\pm~~\mbox{in}~~L^{q}(\mathbb{R}^{N}),~~\forall q\geq N\; \mbox{and}\:
		w_{n}^\pm \rightarrow w^\pm ~~\mbox{a.e. in }~~\mathbb{R}^{N}.
		\end{equation*}
	
Note that, according to \eqref{eq:2.7}, we have
	\begin{equation*}\label{eq:3.21}
	f(x, w_n^\pm)~ w_n^\pm  \leq \e | w_n^\pm|^N + C |w_n^\pm|^{q} exp(\alpha |w_n^\pm|^{\gamma})=:g(w_n^\pm(x)), \quad \mbox{for all}\; \alpha > \alpha_0\, \mbox{and}\, q>N.
	\end{equation*}
	We will prove that $g(w_n^\pm(x))$ is convergent in $L^1(\mathbb{R}^{N})$. First note that
	\begin{equation}\label{eq:3.3''}
	|w_{n}|^N \rightarrow |w|^N~~ \mbox{in}~~ L^{1}(B).
	\end{equation}
	Considering $s,s' > 1$ such that $\frac{1}{s} +\frac{1}{s'}=1$ and $s$ close to $1$, we get
	\begin{equation}\label{eq:4.9}
	|w_{n}|^q \rightarrow |w|^q~~ \mbox{in}~~ L^{s'}(\mathbb{R}^{N}).
	\end{equation}	
	On the other hand, using Lemma \ref{lem8}, we obtain that
	\begin{equation} \label{eq:4.10}
	\begin{aligned}
	\displaystyle c_{\mathcal{N}}&=\limsup_{n \rightarrow +\infty}  \mathcal{E}(w_n)\\
	& = \lim_{n \rightarrow +\infty} \big(\mathcal{E}(w_n)-\frac{1}{\theta}\langle \mathcal{E}'(w_n),w_n\rangle \big) \\
	& =\limsup_{n \rightarrow +\infty}\big(  \frac{1}{N}M(\| w_n\|^{N})-\frac{1}{\theta} m(\| w_n\|^{N})\| w_n\|^{N} +\frac{1}{\theta}\int_{\mathbb{R}^{N}}\big(f(x, w_n) w_n- \theta F(x, w_n)\big)dx\big)  \\
	& > \limsup_{n \rightarrow +\infty} (\frac{1}{2N} - \frac{1}{\theta})m_{0}\|w_{n}\|^{N},
	\end{aligned}
	\end{equation}
	which together with Lemma \ref{lem9} gives that $\displaystyle \limsup_{n \rightarrow +\infty} \|w_n\|^{\gamma} < \frac{\alpha_{N,\beta}}{2(\alpha_0+\delta)}.$
	
	Now choosing $\alpha= \alpha_0+\delta$, we get that
	\begin{equation}\label{eq:4.11}\begin{array}{llll}\displaystyle
	\int_{\mathbb{R}^{N}} \exp(s \alpha |w_n^\pm|^{\gamma})dx &\leq&\displaystyle \int_{\mathbb{R}^{N}}\exp\big(s(\alpha_{0} + \delta )\|w_{n}\|^{\gamma}(\frac{w_{n}}{\|w_{n}\|})^{\gamma}\big)dx\\& \leq&  \displaystyle \int_{\mathbb{R}^{N}}\exp\big(s(\alpha_{0} + \delta )\frac{\alpha_{N,\beta}}{2(\alpha_0+\delta)}(\frac{w_{n}}{\|w_{n}\|})^{\gamma}\big)dx.\end{array}
	\end{equation}
	Since $s > 1$ and is sufficiently close to $1,$ we get $\frac{s}{2}\alpha_{N,\beta}\leq \alpha_{N,\beta}.$ Then it follows by Theorem \ref{th3} that there is $M > 0$ such that
	\begin{equation}\label{eq:4.12}
	\int_{\mathbb{R}^{N}} \exp(s \alpha |w_n^\pm|^{\gamma})dx \leq M.
	\end{equation}
	Since
	\begin{equation}\label{eq:4.13}
	\exp(\alpha |w_n^\pm|^{\gamma}) \rightarrow \exp(\alpha |w^\pm|^{\gamma})\:  \mbox{a.e in }~~\mathbb{R}^{N}.
	\end{equation}
	Then, from \eqref{eq:4.12} and [\cite{Ka}, Lemma 4.8, Chapter 1], we get that
	\begin{equation}\label{eq:4.14}
	\exp(\alpha |w_n^\pm|^{\gamma}) \rightharpoonup \exp(\alpha |w^\pm|^{\gamma})\:  \mbox{ in }~~L^{s}(\mathbb{R}^{N}).
	\end{equation}
	Now using \eqref{eq:3.3''}, \eqref{eq:4.9}, \eqref{eq:4.14} and by proceeding as in the lemma \ref{lemma7}, we will complete the proof of  \eqref{eq:3.71}.\\

$(ii)$	Now, proceeding in the similar way to the proof of lemma \ref{lemma7}, there exists $w \in \mathcal{N}$ such that $\mathcal{E}(w) = c_{\mathcal{N}}$, which is the conclusion we want.	
\end{proof}
Therefore, from Lemma \ref{lemma4}, we deduce that $w$ is a least energy sign-changing solution for problem \eqref{eq:1.1}.
Furthermore, in the case when $m(t)=a+bt$,  $a,b>0$, from Lemma \ref{lemma8}, we deduce that $w$ is a least energy sign-changing solution for problem \eqref{eq:1.1} with exactly two nodal domains.

\subsection{Funding Information }
No funding
\subsection{Author Contribution}
The authors contributed equally to this work.
\subsection{Conflict of interest}
The authors declare that they have no conflicts of interests.


\begin{thebibliography}{99.}
\bibitem{ABJ}{ Abid I., Baraket S. and Jaidane R.,} \emph{On a weighted elliptic equation of N-Kirchhoff type} Accepted in 2022  to appear in Demonstratio Mathematica.

\bibitem{Adi1}{Adimurthi A. and Sandeep K.,} \emph{A Singular Moser-Trudinger Embedding and Its Applications,} Nonlinear Differential Equations and Applications, 2007, vol. 13, issue 5-6, 585-603. DOI:10.1007/s00030-006-4025-9

\bibitem{AC} { Alves C.O. and  Corrêa F.J.S.A.,} \emph{On existence of solutions for a class of problem involving a nonlinear operator,} Comm. Appl. Nonlinear Anal. 8 (2001), 43-56.

 \bibitem{ACM} {Alves C.O., Corrêa F.J.S.A. and  Ma T.F.,} \emph{Positive solutions for a quasilinear elliptic equation of Kirchhoff type,} Comput. Math. Appl. 49 (2005), 85-93.

\bibitem{ABC} { Ambrosetti A.,  Badiale M. and  Cingolani S.,} \emph{Semiclassical stats of nonlinear Schr\"{o}dinger equations with potentials,} Arch. Ration. Mech. Anal. 140, 285–300 (1997).

\bibitem{A M S} {Ambrosetti A.,  Malchiodi A. and  Secchi S.,} \emph{Multiplicity results for some nonlinear Schr\"{o}dinger equations with potentials,} Arch. Ration. Mech. Anal. 159, 253–271 (2001).

\bibitem{AOJ}{Sami Aouaoui , Rahma Jlel,}\emph{New weighted sharp Trudinger–Moser inequalities defined
on the whole euclidean space RN and applications,} Calc. Var. (2021) 60:50
https://doi.org/10.1007/s00526-021-01925-7
\bibitem{BJ} {Baraket S. and  Jaidane R.,} \emph{Non-autonomous weighted elliptic equations with double exponential  growth,,} An. \c{S}t. Univ. Ovidius Constan\c{t}a, Vol. 29(3),2021, 33-66.

\bibitem{CLM1} {Caglioti E., Lions P.L.,  Marchioro C. and Pulvirenti M.,} \emph{A Special Class of Stationary Flows for Two-Dimensional Euler Equations: a Statistical Mechanics Description,} Communications in Mathematical Physics, 1992, vol. 143, no. 3, 501-525.

\bibitem{CLM2} {Caglioti E., Lions P.L.,  Marchioro C. and  Pulvirenti M.,} \emph{A Special Class of Stationary Flows for Two-Dimensional Euler Equations: a Statistical Mechanics Description. II,} Communications in Mathematical Physics, 1995, vol. 174, no. 2, pp. 229-260.

\bibitem{CR1} {Calanchi M. and  Ruf B.,} \emph{On a Trudinger-Moser type inequalities with logarithmic weights}, Journal of Differential Equations no. 3 (2015), 1967-1989.

\bibitem{CR2} {Calanchi M. and Ruf B.,} \emph{Trudinger-Moser type inequalities with logarithmic weights in dimension N,} Nonlinear Analysis, Series A; Theory Methods and Applications 121 (2015), 403-411. DOI: 10.1016/j.na.2015.02.001.

\bibitem{CR3} {Calanchi M. and Ruf B.,} \emph{Weighted Trudinger-Moser inequalities and Applications,} Bulletin of the South Ural State University. Ser. Mathematical Modelling, programming and Computer Software vol. 8  no. 3 (2015), 42-55. DOI: 10.14529/mmp150303

\bibitem{CRS} {Calanchi M.,  Ruf B. and Sani F.,} \emph{Elliptic equations in dimension $2$ with double exponential nonlinearities,} NoDea Nonlinear Differ. Equ. Appl., 24 (2017), Art. 29. DOI: 10.1007/s00030-017-0453-y

\bibitem{CT} {Calanchi M. and Terraneo E.,} \emph{Non-radial Maximizers For Functionals With Exponential Non- linearity in $\mathbb{R}^{2}$,} Advanced Nonlinear Studies vol. 5 (2005), 337-350. DOI:10.1515/ans-2005-0302

\bibitem{CK} { Chanillo S. and  Kiessling M.,} \emph{Rotational Symmetry of Solutions of Some Nonlinear Problems in Statistical Mechanics and in Geometry,} Communications in Mathematical Physics, 1994, vol. 160, no. 2,  217-238. DOI: 10.1007/BF02103274

\bibitem{CTW} {Chen S. T., Tang X. H. and Wei J. Y.,} \emph{Improved results on planar Kirchhoff-type elliptic problems with critical exponential growth,} Z. Angew. Math. Phys. 72, 38 (2021). 

\bibitem{CL}{ Cingolani S. and Lazzo M.,} \emph{Multiple positive solutions to nonlinear Schr\"{o}dinger equations with competing potential functions,} J. Differ. Equations 160, 118–138 (2000).

\bibitem{FMR}{ de Figueiredo D.G.,  Miyagaki O.H. and  Ruf B.,} \emph{Elliptic equations in $\mathbb{R}^{2}$ with nonlinearities in the critical growth range,} Calc. Var. Partial Differential Equations 3 (2) (1995), 139-153. DOI: 10.1007/BF01205003.

\bibitem{FN}{ Figueiredo G. M. and  Nascimento R. G.,} \emph{Existence of a nodal solution with minimal energy for a Kirchhoff equation,} Math. Nachr. 288 (2015), 48-60.

\bibitem{FN1}{Figueiredo G. M. and Nunes F.,} \emph{Existence of positive solutions for a class of quasilinear elliptic problems with exponential growth via the Nehari manifold method,} Rev. Mat. Complut. 32 (2019), 1-18.

\bibitem{FJR}{de Figueiredo D.G., do Ó J.M. and  Ruf B.,} \emph{On an inequality by N. Trudinger and J. Moser and related elliptic equations,} Comm. Pure Appl. Math. LV (2002) 135-152.

\bibitem{ddr} { de Figueiredo D.G.,  do  \'O J.M. and  Ruf  B.,} \emph{Elliptic equations and systems with critical Trudinger-Moser nonlinearities,} Discrete Contin. Dyn. Syst., 30(2), (2011) 455-476.

\bibitem{FSJ}{Figueiredo G. M. and Santos Junior J. R.,} \emph{Existence of a least energy nodal solution for a Schr\"{o}dinger-Kirchhoff equation with potential vanishing at infinity,} J. Math. Phys. 56 (2015), 051506.

\bibitem{FS}{ de Figueiredo D. G. and Uberlandio B. Severo,} \emph{Ground State Solution for a Kirchhoff Problem with Exponential Critical Growth,} Milan J. Math. (2015)DOI 10.1007/s00032-015-0248-8

\bibitem{DKN} {Drabek P.,  Kufner A. and Nicolosi F.,} \emph{Quasilinear Elliptic Equations with Degenerations and Singularities}, Walter de Gruyter, Berlin (1997). DOI:10.1515/9783110804775

\bibitem{DJ1} {Dridi B. and  Jaidane R.,} \emph{Nodal solutions for logarithmic weighted N-laplacian problem with exponential nonlinearities,} Ann Univ Ferrara (2023). https://doi.org/10.1007/s11565-023-00457-6.

\bibitem{DJ2} {Dridi B. and  Jaidane R.,} \emph{ Existence Solutions for a Weighted Biharmonic Equation with Critical Exponential Growth.} Mediterr. J. Math. 20, 102 (2023). https://doi.org/10.1007/s00009-023-02301-9

\bibitem{DV}{ Dumitru Motreanu N. P. and  Viorica Venera Motreanu,} \emph{Topological and variational methods with applications to nonlinear boundary value problems.} Springer, New York, NY, 2014. 13

\bibitem{GCZ}{Gao L., Chen C.F. and Zhu C.X.,} \emph{Existence of sign-changing solutions for Kirchhoff equations with critical or supercritical nonlinearity,} Appl. Math. Lett. 107 (2020), Article 106424.

\bibitem{HY} { Han W. and Yao J.,}  \emph{The sign-changing solutions for a class of p-Laplacian Kirchhoff type problem in bounded domains,} Comput. Math. Appl. 76 (2018), 1779–1790.

\bibitem{KJ}{Kharrati S. and Jaidane R.,} \emph{ Existence of Positive Solutions to Weighted Linear Elliptic Equations Under Double Exponential Nonlinearity Growth.} Bull. Iran. Math. Soc. 48, 993–1021 (2022). https://doi.org/10.1007/s41980-021-00559-x

\bibitem{Ka}{Kavian O.,} \emph{ Introduction à la Théorie des Points Critiques. Springer-Verlag,} Berlin, 1991. 15, 17

\bibitem{K} { Kiessling M.K.-H.,} \emph{Statistical Mechanics of Classical Particles with Logarithmic Interactions,} Communications on Pure and Applied Mathematics, 1993, vol. 46, 27-56. DOI:10.1002/cpa.3160460103

\bibitem{ki} {Kirchhof G.,} \emph{Mechanik}, Teubner, Leipzig, 1883.

\bibitem{Kuf} {Kufner A.,} \emph{Weighted Sobolev spaces,} John Wiley and Sons Ltd, 1985.  Doi: 10.1112/blms/18.2.220

\bibitem{LDZ}{ Li Q., Du X. and Zhao  Z.,} \emph{Existence of sign-changing solutions for nonlocal Kirchhoff-Schrödinger-type equations in $\R^3$,} J. Math. Anal. Appl. 477 (2019), 174–186.

\bibitem{liang4}{Liang S. and R\u{a}dulescu V. D.,} \emph{Least-energy nodal solutions of critical Kirchhoff problems with logarithmic nonlinearity}, Anal. Math. Phys., {10:45} (2020), 1--31.

\bibitem{Li}{Lions J.,}\emph{On some questions in boundary value problems of mathematical physics,} North-Holland Math. Stud. 30, 284–346 (1978).

 \bibitem{Liou} {Liouville J.,} \emph{Sur l'equation aux deriv\'{e}es partielles,} Journal de Math\'{e}matiques Pures et Appliqu\'{e}es, 1853, vol. 18, 71-72.

 \bibitem{Lu}{Lu G. and Yang Y.,} \emph{Sharp constant and extremal function for the improved Moser-Trudinger inequality involving lp norm in two dimension,}  Discrete Contin. Dyn. Syst., 25(3), (2012) 963-979.

\bibitem{MS} {Masmoudi N. and Sani F.,} \emph{Trudinger-Moser inequalities with the exact growth condition in $\mathbb{R^{N}}$ and applications,} Comm. Partial Diferential Equations 40 no. 8 (2015), 1408-1440. DOI:10.1080/03605302.2015.1026775

\bibitem{mi}{Miranda C.,} \emph{Un'osservazione su un teorema di Brouwer}, Boll Un Mat Ital, 3 (1940), 5-7.

 \bibitem{Mos}{Moser J.,}  \emph{A sharp form of an inequality by N. Trudinger,} Indiana Univ. Math. J., 20(71), (1970) 1077-1092.

\bibitem{Po} {Pohozaev S.I.,} \emph{The Sobolev Embedding in the Case $pl = n$.} Proc. of the Technical Scientific Conference on Advances of Scientific Research, 1964-1965, Mathematics Section, (Moskov. Energet. Inst., Moscow), pp. 158-170.

\bibitem{S} {Shen L.,} \emph{Sign-changing solutions to a N-Kirchhoff equation with critical exponential growth in $\R^N$,} Bull. Malays. Math. Sci. Soc., 44 (2021), 3553–3570.

\bibitem{Sh} {Shuai W.,} \emph{ Sign-changing solutions for a class of Kirchhoff-type problem in bounded domains,} J. Differ. Equ. 259 (2015), 1256–1274.


\bibitem{SSLR} {de Souza M., Severo U. B. and Luiz do Rego T.,} \emph{On solutions for a class of fractional Kirchhoff-type problems with Trudinger–Moser nonlinearity,} Comm. Contemp. Math.,Vol. 24, No. 05, 2150002 (2022).


\bibitem{T1} { Tarantello G.,} \emph{Condensate Solutions for the Chern - Simons - Higgs Theory,} Journal of Mathematical Physics vol. 37  (1996), 3769-3796. DOI: 10.1063/1.531601

\bibitem{T2} {Tarantello G.,} \emph{Analytical Aspects of Liouville-Type Equations with Singular Sources,} Handbook of Differential Equations (M. Chipot and P. Quittner, eds.), Elsevier, North Holland 2004, 491-592.

\bibitem{Tru} {Trudinger N.,} \emph{On imbeddings into Orlicz spaces and some applications,} J. Math. Mech., 17, (1967) 473-483.


\bibitem{WTC} { Wen L., Tang X. H. and Chen S.,} \emph{Ground state sign-changing solutions for Kirchhoff equations with
logarithmic nonlinearity,} Electron. J. Qual. Theor., 47 (2019), 1–13.

\bibitem{mw1996} { Willem M.,} \emph{Minimax Theorems}, Birkh\"auser, Boston, 1996.

\bibitem{XTZ} {Xiao T., Tang Y. and Zhang Q.,} \emph{The existence of sign-changing solutions for Schr\"{o}dinger-Kirchhoff problems in $\R^3$,} AIMS Mathematics, 6(7), (2021), 6726-6733.


\bibitem{Ya} { Yudovich V.I.,} \emph{Some Estimates Connected with Integral Operators and with Solutions of Elliptic Equations.} Dokl. Akad. Nauk SSSR, 138, (1961) pp. 805-808.

\bibitem{ZYL} {Zhang Y., Yang Y. and Liang S.,} \emph{Least energy sign-changing solution for N-Laplacian problem with logarithmic and exponential nonlinearities,} Journal of Mathematical Analysis and Applications. 505(1): 125432, 2022.




























\end{thebibliography}
\end{document}